\documentclass[12pt]{article}
\usepackage{latexsym, amssymb}
\usepackage{mathbbol}

\usepackage{typearea}
\typearea{15}
\usepackage{graphicx}
\usepackage{tikz}
\usepackage{ifthen}
\usepackage{caption}

\usetikzlibrary{positioning,shapes,calc}

\DeclareMathAlphabet\gothic{U}{euf}{m}{n}

\makeatletter
\def\eqnarray{\stepcounter{equation}\let\@currentlabel=\theequation
\global\@eqnswtrue
\tabskip\@centering\let\\=\@eqncr
$$\halign to \displaywidth\bgroup\hfil\global\@eqcnt\z@
  $\displaystyle\tabskip\z@{##}$&\global\@eqcnt\@ne
  \hfil$\displaystyle{{}##{}}$\hfil
  &\global\@eqcnt\tw@ $\displaystyle{##}$\hfil
  \tabskip\@centering&\llap{##}\tabskip\z@\cr}

\def\endeqnarray{\@@eqncr\egroup
      \global\advance\c@equation\m@ne$$\global\@ignoretrue}

\def\@yeqncr{\@ifnextchar [{\@xeqncr}{\@xeqncr[5pt]}}
\makeatother

\begin{document}
\bibliographystyle{tom}

\newtheorem{lemma}{Lemma}[section]
\newtheorem{thm}[lemma]{Theorem}
\newtheorem{cor}[lemma]{Corollary}
\newtheorem{voorb}[lemma]{Example}
\newtheorem{rem}[lemma]{Remark}
\newtheorem{prop}[lemma]{Proposition}
\newtheorem{ddefinition}[lemma]{Definition}
\newtheorem{stat}[lemma]{{\hspace{-5pt}}}

\newenvironment{remarkn}{\begin{rem} \rm}{\end{rem}}
\newenvironment{exam}{\begin{voorb} \rm}{\end{voorb}}
\newenvironment{definition}{\begin{ddefinition} \rm}{\end{ddefinition}}

\newcommand{\gota}{\gothic{a}}
\newcommand{\gotb}{\gothic{b}}
\newcommand{\gotc}{\gothic{c}}
\newcommand{\gote}{\gothic{e}}
\newcommand{\gotf}{\gothic{f}}
\newcommand{\gotg}{\gothic{g}}
\newcommand{\gothh}{\gothic{h}}
\newcommand{\gotk}{\gothic{k}}
\newcommand{\gotm}{\gothic{m}}
\newcommand{\gotn}{\gothic{n}}
\newcommand{\gotp}{\gothic{p}}
\newcommand{\gotq}{\gothic{q}}
\newcommand{\gotr}{\gothic{r}}
\newcommand{\gots}{\gothic{s}}
\newcommand{\gotu}{\gothic{u}}
\newcommand{\gotv}{\gothic{v}}
\newcommand{\gotw}{\gothic{w}}
\newcommand{\gotz}{\gothic{z}}
\newcommand{\gotA}{\gothic{A}}
\newcommand{\gotB}{\gothic{B}}
\newcommand{\gotG}{\gothic{G}}
\newcommand{\gotL}{\gothic{L}}
\newcommand{\gotS}{\gothic{S}}
\newcommand{\gotT}{\gothic{T}}

\newcounter{teller}
\renewcommand{\theteller}{(\alph{teller})}
\newenvironment{tabel}{\begin{list}%
{\rm  (\alph{teller})\hfill}{\usecounter{teller} \leftmargin=1.1cm
\labelwidth=1.1cm \labelsep=0cm \parsep=0cm}
                      }{\end{list}}

\newcounter{tellerr}
\renewcommand{\thetellerr}{{\upshape{(\roman{tellerr})}}}
\newenvironment{tabeleq}{\begin{list}%
{\rm  (\roman{tellerr})\hfill}{\usecounter{tellerr} \leftmargin=1.1cm
\labelwidth=1.1cm \labelsep=0cm \parsep=0cm}
                         }{\end{list}}

\newcounter{tellerrr}
\renewcommand{\thetellerrr}{(\Roman{tellerrr})}
\newenvironment{tabelR}{\begin{list}%
{\rm  (\Roman{tellerrr})\hfill}{\usecounter{tellerrr} \leftmargin=1.1cm
\labelwidth=1.1cm \labelsep=0cm \parsep=0cm}
                         }{\end{list}}

\newcounter{proofstep}
\newcommand{\nextstep}{\refstepcounter{proofstep}\ruimte \par 
          \noindent{\bf Step \theproofstep} \hspace{5pt}}
\newcommand{\firststep}{\setcounter{proofstep}{0}\nextstep}

\newcommand{\Ni}{{\mathbb{N}}}
\newcommand{\Qi}{{\mathbb{Q}}}
\newcommand{\Ri}{{\mathbb{R}}}
\newcommand{\Ci}{{\mathbb{C}}}
\newcommand{\Ti}{{\mathbb{T}}}
\newcommand{\Zi}{{\mathbb{Z}}}
\newcommand{\Fi}{{\mathbb{F}}}

\newcommand{\proof}{\mbox{\bf Proof} \hspace{5pt}} 
\newcommand{\remark}{\mbox{\bf Remark} \hspace{5pt}}
\newcommand{\ruimte}{\vskip10.0pt plus 4.0pt minus 6.0pt}

\newcommand{\simh}{{\stackrel{{\rm cap}}{\sim}}}
\newcommand{\ad}{{\mathop{\rm ad}}}
\newcommand{\Ad}{{\mathop{\rm Ad}}}
\newcommand{\Aut}{\mathop{\rm Aut}}
\newcommand{\arccot}{\mathop{\rm arccot}}
\newcommand{\capp}{{\mathop{\rm cap}}}
\newcommand{\rcapp}{{\mathop{\rm rcap}}}
\newcommand{\diam}{\mathop{\rm diam}}
\newcommand{\divv}{\mathop{\rm div}}
\newcommand{\ess}{{\rm ess}}
\newcommand{\codim}{\mathop{\rm codim}}
\newcommand{\RRe}{\mathop{\rm Re}}
\newcommand{\IIm}{\mathop{\rm Im}}
\newcommand{\Tr}{{\mathop{\rm Tr \,}}}
\newcommand{\Trind}[1]{{{\mathop{\rm Tr}}_{#1} \,}}
\newcommand{\Vol}{{\mathop{\rm Vol}}}
\newcommand{\card}{{\mathop{\rm card}}}
\newcommand{\supp}{\mathop{\rm supp}}
\newcommand{\sgn}{\mathop{\rm sgn}}
\newcommand{\essinf}{\mathop{\rm ess\,inf}}
\newcommand{\esssup}{\mathop{\rm ess\,sup}}
\newcommand{\Int}{\mathop{\rm Int}}
\newcommand{\lcm}{\mathop{\rm lcm}}
\newcommand{\loc}{{\rm loc}}

\newcommand{\mod}{\mathop{\rm mod}}
\newcommand{\spann}{\mathop{\rm span}}
\newcommand{\one}{\mathbb{1}}

\hyphenation{groups}
\hyphenation{unitary}

\newcommand{\tfrac}[2]{{\textstyle \frac{#1}{#2}}}

\newcommand{\ca}{{\cal A}}
\newcommand{\cb}{{\cal B}}
\newcommand{\cc}{{\cal C}}
\newcommand{\cd}{{\cal D}}
\newcommand{\ce}{{\cal E}}
\newcommand{\cf}{{\cal F}}
\newcommand{\ch}{{\cal H}}
\newcommand{\ci}{{\cal I}}
\newcommand{\ck}{{\cal K}}
\newcommand{\cl}{{\cal L}}
\newcommand{\cm}{{\cal M}}
\newcommand{\co}{{\cal O}}
\newcommand{\cs}{{\cal S}}
\newcommand{\ct}{{\cal T}}
\newcommand{\cx}{{\cal X}}
\newcommand{\cy}{{\cal Y}}
\newcommand{\cz}{{\cal Z}}

\newlength{\hightcharacter}
\newlength{\widthcharacter}
\newcommand{\covsup}[1]{\settowidth{\widthcharacter}{$#1$}\addtolength{\widthcharacter}{-0.15em}\settoheight{\hightcharacter}{$#1$}\addtolength{\hightcharacter}{0.1ex}#1\raisebox{\hightcharacter}[0pt][0pt]{\makebox[0pt]{\hspace{-\widthcharacter}$\scriptstyle\circ$}}}
\newcommand{\cov}[1]{\settowidth{\widthcharacter}{$#1$}\addtolength{\widthcharacter}{-0.15em}\settoheight{\hightcharacter}{$#1$}\addtolength{\hightcharacter}{0.1ex}#1\raisebox{\hightcharacter}{\makebox[0pt]{\hspace{-\widthcharacter}$\scriptstyle\circ$}}}
\newcommand{\scov}[1]{\settowidth{\widthcharacter}{$#1$}\addtolength{\widthcharacter}{-0.15em}\settoheight{\hightcharacter}{$#1$}\addtolength{\hightcharacter}{0.1ex}#1\raisebox{0.7\hightcharacter}{\makebox[0pt]{\hspace{-\widthcharacter}$\scriptstyle\circ$}}}
\newcommand{\excov}[1]{\settowidth{\widthcharacter}{$#1$}\addtolength{\widthcharacter}{-0.15em}\settoheight{\hightcharacter}{$#1$}\addtolength{\hightcharacter}{4mm}#1\raisebox{\hightcharacter}{\makebox[0pt]{\hspace{-\widthcharacter}$\scriptstyle\circ$}}}

\thispagestyle{empty}

\vspace*{1cm}
\begin{center}
{\Large\bf The Dirichlet-to-Neumann operator on rough domains} \\[5mm]

\large W. Arendt$^1$ and A.F.M. ter Elst$^2$

\end{center}

\vspace{5mm}

\begin{center}
{\bf Abstract}
\end{center}

\begin{list}{}{\leftmargin=1.8cm \rightmargin=1.8cm \listparindent=10mm 
   \parsep=0pt}
\item

We consider a bounded connected open set $\Omega \subset \Ri^d$ whose boundary 
$\Gamma$ has a finite $(d-1)$-dimensional Hausdorff measure.
Then we define the Dirichlet-to-Neumann operator $D_0$ on 
$L_2(\Gamma)$ by form methods.
The operator $-D_0$ is self-adjoint and generates a contractive 
$C_0$-semigroup $S = (S_t)_{t > 0}$ on $L_2(\Gamma)$.
We show that the asymptotic behaviour of $S_t$ as $t \to \infty$
is related to properties of the trace of functions in $H^1(\Omega)$
which $\Omega$ may or may not have.

\end{list}

\vspace{4cm}
\noindent
April 2010

\vspace{5mm}

\vspace{5mm}

\vspace{15mm}

\noindent
{\bf Home institutions:}    \\[3mm]
\begin{tabular}{@{}cl@{\hspace{10mm}}cl}
1. & Institute of Applied Analysis  & 
  2. & Department of Mathematics   \\
& University of Ulm   & 
  & University of Auckland   \\
& Helmholtzstr.\ 18 & 
  & Private bag 92019  \\
& 89081 Ulm & 
  & Auckland  \\
& Germany  & 
  & New Zealand  \\[8mm]
\end{tabular}


\mbox{}
\thispagestyle{empty}

\newpage
\setcounter{page}{1}

\section{Introduction} \label{Sdton1}

Throughout this paper $\Omega$ is a bounded, connected, open set in 
$\Ri^d$ with boundary~$\Gamma$.
We consider the $(d-1)$-dimensional Hausdorff measure $\sigma$ on $\Gamma$,
where $d \geq 2$
and assume throughout that $\sigma(\Gamma) < \infty$.
The purpose of this article is to define the Dirichlet-to-Neumann
operator $D_0$ on $L_2(\Gamma)$ under these mild assumptions on $\Omega$
and to study the semigroup $(S_t)_{t > 0}$ generated by $-D_0$ 
on $L_2(\Gamma)$.

For this purpose we define at first the trace in the following way.
Given $u \in H^1(\Omega)$, a function $\varphi \in L_2(\Gamma)$ is 
called a {\bf trace} of $u$ if there exists a sequence $(u_n)_{n \in \Ni}$
in $H^1(\Omega) \cap C(\overline \Omega)$ such that 
$\lim_{n \to \infty} u_n = u$ in $H^1(\Omega)$ and 
$\lim_{n \to \infty} u_n|_\Gamma = \varphi$ in $L_2(\Gamma)$.
The trace may not be unique (see Example~\ref{xdton450}).
If $u$ has a trace, then $u \in \widetilde H^1(\Omega)$, the 
closure of $H^1(\Omega) \cap C(\overline \Omega)$ in $H^1(\Omega)$.
The space $\widetilde H^1(\Omega)$ may be a proper subset of $H^1(\Omega)$
and in general not every $u \in \widetilde H^1(\Omega)$ has a trace.

Next we define the (weak) normal derivative via Green's formula.
Let $u \in H^1(\Omega)$ be such that $\Delta u \in L_2(\Omega)$
as distribution.
We say that $u$ has a {\bf normal derivative in $L_2(\Gamma)$} if 
there exists a $\psi \in L_2(\Gamma)$ such that 
\[
\int_\Omega (\Delta u) \, v  + \int_{\Omega} \nabla u \cdot \nabla v
 = \int_\Gamma \psi \, v \, d\sigma
\]
for all $v \in H^1(\Omega) \cap C(\overline \Omega)$.
In that case $\psi$ is unique. 
We set 
$\frac{\partial u}{\partial \nu} := \psi$ and call 
it the normal derivative of $u$.
Now we define the Dirichlet-to-Neumann operator $D_0$ on $L_2(\Gamma)$
as follows.
Given $\varphi,\psi \in L_2(\Gamma)$, we say that $\varphi \in D(D_0)$ and 
$D_0 \varphi = \psi$ if there exists a $u \in H^1(\Omega)$ such that 
$\Delta u = 0$ as distribution, 
$\varphi$ is a trace of $u$,
the function $u$ has a normal derivative in $L_2(\Gamma)$
and $\frac{\partial u}{\partial \nu} = \psi$.
Even though the function $u$ might not have a unique trace,
we shall prove that the operator $D_0$ is univocal.
In fact, $D_0$ is a self-adjoint operator on $L_2(\Gamma)$ and 
$-D_0$ generates a positive $C_0$-semigroup $S$ on $L_2(\Gamma)$
satisfying $S_t \one_\Gamma = \one_\Gamma$ for all $t > 0$.
This is true without any regularity hypothesis on $\Omega$
(besides $\sigma(\Gamma) < \infty$).
One purpose of this paper is to show that diverse properties concerning the asymptotic 
behaviour of $S_t$ as $t \to \infty$ are related to properties of the
trace, which in fact are properties of $\Omega$, which $\Omega$ 
may or may not have.

Here are our main results.

\paragraph*{A. Strong convergence of $S$. }
We say that {\bf the trace on $\Omega$ is unique} if the function 
$\varphi = 0 \in L_2(\Gamma)$ is the only trace of $u = 0 \in H^1(\Omega)$.
Then every element of $\widetilde H^1(\Omega)$ has at most one trace.
This is true in many cases, but 
not in general.
Define $P \colon L_2(\Gamma) \to L_2(\Gamma)$ by 
$P \varphi = \Big( \frac{1}{\sigma(\Gamma)} \int_\Gamma \varphi  \Big) \one_\Gamma$.

\begin{thm} \label{tdton101}
The following are equivalent.
\begin{tabeleq}
\item \label{tdton101-1}
The trace on $\Omega$ is unique.
\item \label{tdton101-2}
$\dim(\ker D_0) = 1$.
\item \label{tdton101-3}
$\lim_{t \to \infty} S_t \varphi = P \varphi$ for all 
$\varphi \in L_2(\Gamma)$.
\item \label{tdton101-4}
$S$ is irreducible.
\end{tabeleq}
\end{thm}

The irreducibility of $S$ is surprising since the boundary $\Gamma$
need not be connected (consider an annulus for example).
Thus this result reflects somehow that the operator $D_0$ is not local.

\paragraph*{B. Norm convergence of $S$. }
We emphasize that in general not every element in $\widetilde H^1(\Omega)$
has a trace and if it has a trace, then it might not be unique.
We next characterize when both properties are valid, i.e.\ every element of $\widetilde H^1(\Omega)$
has a trace and this trace is unique.
This is true for example if $\Omega$ has a Lipschitz boundary.

\begin{thm} \label{tdton102}
The following are equivalent.
\begin{tabeleq}
\item \label{tdton102-1}
$\lim_{t \to \infty} S_t = P$ in $\cl(L_2(\Gamma))$.
\item \label{tdton102-2.5}
There exists a $c > 0$ such that 
\[
\int_\Gamma |u|^2 \leq c \int_\Omega |\nabla u|^2
\]
for all $u \in H^1(\Omega) \cap C(\overline \Omega)$ with 
$\int_\Gamma u = 0$.
\item \label{tdton102-2}
There exists a $c > 0$ such that 
\[
\int_\Gamma |u|^2 
\leq c \Big( \int_\Omega |\nabla u|^2 + \int_\Omega |u|^2 \Big)
\]
for all $u \in H^1(\Omega) \cap C(\overline \Omega)$.
\item \label{tdton102-3}
Every $u \in \widetilde H^1(\Omega)$ has a unique trace.
\item \label{tdton102-4}
$0 \not\in \sigma_\ess(D_0)$.
\end{tabeleq}
\end{thm}

\paragraph*{C. Compactness of the resolvent. }
We shall show that the operator $D_0$ has compact resolvent if and only if every 
$u \in \widetilde H^1(\Omega)$ has a unique trace $\Tr u$ and 
the map $\Tr \colon \widetilde H^1(\Omega) \to L_2(\Gamma)$ is 
compact.
This implies that the embedding $\widetilde H^1(\Omega) \to L_2(\Omega)$
is also compact.
We construct, however, a bounded domain with continuous boundary and 
with $\sigma(\Gamma) < \infty$, such that $D_0$ does not have compact 
resolvent (even though the embedding 
$H^1(\Omega) = \widetilde H^1(\Omega) \hookrightarrow L_2(\Omega)$ is 
compact since the boundary is continuous).

\ruimte

The Dirichlet-to-Neumann operator is a well-known object occurring in 
many applications. 
In general it is considered on domains of class $C^\infty$, though,
see e.g.\ the monograph of Taylor \cite{Tay3} \cite{Tay5} \cite{Tay6}.
Then the operator fits into the framework of pseudo-differential operators
and also semigroup properties are studied \cite{Esc1} \cite{Eng}.
Our point is the very general variational definition which allows an 
easy approach also for rough domains.
On the other hand, the questions concerning trace properties which we investigate
here become delicate.
They are the main subject of the paper.
Some of the trace properties considered here are related to 
investigations of the Laplace operator with Robin boundary conditions
on arbitrary domains as in \cite{Daners2}, see also \cite{AW2}.

The paper is organized as follows.
In Section~\ref{Sdton1.5} we consider the asymptotic behaviour of Markovian semigroups.
This section is independent of the Dirichlet-to-Neumann operator.
In Section~\ref{Sdton2.5} we prove the existence and uniqueness 
of the Dirichlet-to-Neumann operator on rough domains and 
show that it is a self-adjoint operator which generates a 
Markovian semigroup.
In Section~\ref{Sdton3} we prove Theorem~\ref{tdton101}.
In addition we give other characterizations 
of the uniqueness of the trace in terms of the form 
associated to the Laplacian with Robin boundary conditions and in terms
of the relative capacity.
In Section~\ref{Sdton4} we define the trace as a mapping 
and study its properties.
In Section~\ref{Sdton6} we characterize when every 
element of $\widetilde H^1(\Omega)$ has a trace.
Moreover, we prove Theorem~\ref{tdton102}.
In Section~\ref{Sdton7} we characterize when the 
map $u \mapsto u|_\Gamma$ from 
$(H^1(\Omega) \cap C(\overline \Omega), \|\cdot\|_{H^1(\Omega)})$
into $L_2(\Gamma)$ is compact.
Theorem~\ref{tdton102} and the compactness of the trace 
can be reformulated in terms of the form 
associated to the Laplacian with Robin boundary conditions.
This is done in Section~\ref{Sdton8}.
Finally, in Section~\ref{Sdton5} we present two striking 
examples.

Throughout this paper the field is $\Ri$ and we only consider univocal operators.

\section{Asymptotic behaviour of markovian semigroups} \label{Sdton1.5}

In this section we put together some asymptotic properties of markovian semigroups.
At first we consider a {\bf self-adjoint semigroup}, i.e.\ a semigroup 
consisting of self-adjoint operators.

\begin{prop} \label{pdton1501}
Let $S$ be a contractive $C_0$-semigroup of self-adjoint operators on a Hilbert
space~$H$.
Then 
\[
P_S f = \lim_{t \to \infty} S_t f
\]
exists for all $f \in H$ and $P_S$ is the orthogonal projection onto 
$\ker A$, where $- A$ denotes the generator of $S$.
\end{prop}
\proof\
By the spectral theorem we may assume that $H = L_2(Y)$,
$D(A) = \{ f \in L_2(Y) : m \, f \in L_2(Y) \} $
and $A f = m \, f$ for all $f \in D(A)$, 
where $(Y,\Sigma,\mu)$ is a locally finite measure space
and $m \colon Y \to [0,\infty)$ is a measurable function.
Then $\ker A = \{ f \in L_2(Y) : f = 0 \mbox{ a.e.\ on } Y \setminus Y_0 \} $,
where $Y_0 = m^{-1}( \{ 0 \} )$.
The orthogonal projection $P_S$ onto $\ker A$ is given by 
$P_S f = \one_{Y_0} \, f$.
Moreover, $S_t f = e^{-t \, m} f$ for all $t > 0$ and $f \in L_2(Y)$.
Now the claim follows from the Lebesgue dominated convergence theorem.\hfill$\Box$

\ruimte

Next we consider a finite measure space $(\Gamma,\Sigma,\sigma)$.
A {\bf Markov operator} $T$ on $L_2(\Gamma)$ is an operator 
satisfying $T \one_\Gamma = \one_\Gamma$ and $T f \geq 0$ for all 
$f \in L_2(\Gamma)$ with $f \geq 0$.
As a consequence $T L_\infty(\Gamma) \subset L_\infty(\Gamma)$ 
and $T^{(\infty)} := T|_{L_\infty(\Gamma)}$ is contractive.
If $T$ is a self-adjoint Markov operator on $L_2(\Gamma)$,
then $T$ is contractive for the $L_1$-norm.
Hence for all $p \in [1,\infty]$ there exists a unique 
$T^{(p)} \in \cl(L_p(Y))$ such that 
$T^{(p)} f = T f$ for all $f \in L_p(Y) \cap L_2(Y)$.
Moreover, $\|T^{(p)}\|_{\cl(L_p(Y))} \leq 1$.
The operator $T^{(\infty)}$ is the adjoint of the operator $T^{(1)}$.

A $C_0$-semigroup $S$ on $L_2(\Gamma)$ is called {\bf irreducible}
if for each $\Gamma_1 \in \Sigma$ with 
\[
S_t L_2(\Gamma_1)\subset L_2(\Gamma_1)
\]
for all $t > 0$ it follows that $\sigma(\Gamma_1)=0$ or $\sigma(\Gamma \setminus \Gamma_1)=0$.
Here, and in the sequel, we let
$L_2(\Gamma_1) = \{ f \in L_2(Y) : f = 0 \mbox{ a.e.\ on } \Gamma \setminus \Gamma_1 \} $.
A {\bf Markov semigroup} on $L_2(\Gamma)$ is a $C_0$-semigroup $S$ on $L_2(\Gamma)$ 
such that $S_t$ is a Markov operator for all $t > 0$.
In that case $(S^{(p)}_t)_{t > 0}$ is a positive contractive $C_0$-semigroup 
on $L_p(\Gamma)$ for all $p \in [1,\infty)$.
Moreover, $\Ri \, \one_\Gamma \subset \ker A$, where $-A$ is the generator of $S$.

\begin{prop} \label{pdton1502}
Let $S$ be a self-adjoint Markov semigroup on $L_2(\Gamma)$.
Then $S$ is irreducible if and only if $\ker A = \Ri \, \one_\Gamma$, 
where $-A$ is the generator of $S$.
\end{prop}
\proof\
`$\Rightarrow$'.
This follows from \cite{Nag} Section~C-III, Proposition~3.5(c).

`$\Leftarrow$'.
Let $\Gamma_1 \in \Sigma$ be such that
$S_t L_2(\Gamma_1) \subset L_2(\Gamma_1)$ for all $t>0$.
Set $\Gamma_2:=\Gamma \setminus \Gamma_1$.
Then $L_2(\Gamma_2)=L_2(\Gamma_1)^\perp$ and since $S_t$
is self-adjoint, it follows that $S_t L_2(\Gamma_2) \subset L_2(\Gamma_2)$ for all
$t>0$.
Now 
$\one_{\Gamma_1} + \one_{\Gamma_2}
= \one_\Gamma
= S_t \one_\Gamma
= S_t \one_{\Gamma_1} + S_t \one_{\Gamma_2}$
by assumption.
Moreover, $S_t \one_{\Gamma_j} \in L_2(\Gamma_j)$ vanishes 
outside $\Gamma_j$ for all $j \in \{ 1,2 \} $.
Hence $S_t \one_{\Gamma_1} = \one_{\Gamma_1}$ for all $t > 0$.
This implies that $\one_{\Gamma_1} \in \ker A$.
Since $\ker A = \Ri \, \one_\Gamma$ by assumption,
it follows that $\sigma(\Gamma_1)=0$ or $\sigma(\Gamma_2)=0$.\hfill$\Box$

\ruimte

Next we show that a self-adjoint Markov semigroup is irreducible if and 
only if it converges to an equilibrium.
For all $f \in L_1(\Gamma)$ define 
\begin{equation}
P f = \frac{1}{\sigma(\Gamma)} \Big( \int_\Gamma f \Big) \one_\Gamma
.
\label{eSdton1.5;1}
\end{equation}
Then $P$ defines a positive contractive projection on $L_p(\Gamma)$
for all $p \in [1,\infty]$.

\begin{thm} \label{tdton1.503}
Let $S$ be a self-adjoint Markov semigroup on $L_2(\Gamma)$.
The following are equivalent.
\begin{tabeleq}
\item \label{tdton1.503-1}
$S$ is irreducible.
\item \label{tdton1.503-2}
There exists a $p \in [1,\infty)$ such that 
$\lim_{t \to \infty} S^{(p)}_t f = P f$ in $L_p(\Gamma)$ for all $f \in L_p(\Gamma)$.
\item \label{tdton1.503-3}
For all $p \in [1,\infty)$ one has 
$\lim_{t \to \infty} S^{(p)}_t f = P f$ in $L_p(\Gamma)$ for all $f \in L_p(\Gamma)$.
\end{tabeleq}
\end{thm}
\proof\
`\ref{tdton1.503-1}$\Rightarrow$\ref{tdton1.503-2}'.
If $S$ is irreducible, then $\ker A = \Ri \, \one_\Gamma$ by Proposition~\ref{pdton1502},
where $-A$ is the generator of $S$.
Then the operator $P$ defined in (\ref{eSdton1.5;1}) is the 
orthogonal projection onto $\ker A$.
Then Statement~\ref{tdton1.503-2} follows from Proposition~\ref{pdton1501}.

`\ref{tdton1.503-2}$\Rightarrow$\ref{tdton1.503-3}'.
Let $p \in [1,\infty)$ and suppose that 
$\lim_{t \to \infty} S^{(p)}_t f = P f$ in $L_p(\Gamma)$ for all $f \in L_p(\Gamma)$.
If $f \in L_p(\Gamma)$ then 
$\|S^{(1)}_t f - P f\|_1
\leq (\sigma(\Gamma))^{\frac{1}{p} - 1} \|S^{(p)}_t f - P f\|_p$
for all $t > 0$.
Therefore $\lim_{t \to \infty} S^{(1)}_t f = P f$ in $L_1(\Gamma)$.
Since $L_p(\Gamma)$ is dense in $L_1(\Gamma)$ and 
$ \{ P \} \cup \{ S^{(1)}_t : t > 0 \} $ are uniformly bounded 
in $\cl(L_1(\Gamma))$ it follows that 
$\lim_{t \to \infty} S^{(1)}_t f = P f$ in $L_1(\Gamma)$ for all $f \in L_1(\Gamma)$.

Finally, let $q \in (1,\infty)$. 
If $f \in L_\infty(\Gamma)$ then by interpolation
\[
\|S^{(q)}_t f - P f\|_q
\leq \|S^{(1)}_t f - P f\|_1^\theta \, \|S^{(\infty)}_t f - P f\|_\infty^{1 - \theta}
\leq \|S^{(1)}_t f - P f\|_1^\theta \, (2 \|f\|_\infty)^{1 - \theta}
,  \]
where $\theta = \frac{1}{q}$.
So $\lim_{t \to \infty} S^{(q)}_t f = P f$ in $L_q(\Gamma)$.
Since $L_\infty(\Gamma)$ is dense in $L_q(\Gamma)$ the 
claim follows as before.

`\ref{tdton1.503-3}$\Rightarrow$\ref{tdton1.503-1}'.
Let $f \in \ker A$. 
Then $S_t f = f$ for all $t > 0$.
Consequently $f = Pf \in \Ri \, \one_\Gamma$.
We have shown that $\ker A = \Ri \, \one_\Gamma$.
It follows from Proposition~\ref{pdton1502} that $S$ is irreducible
and \ref{tdton1.503-1} is valid.\hfill$\Box$

\ruimte

If $A$ is a self-adjoint operator, then
$0 \not\in \sigma_\ess(A)$ means by definition that 
$0$ is not an accumulation point of $\sigma(A)$ and $\ker A$ is finite 
dimensional.
Thus if $S$ is a self-adjoint irreducible Markov semigroup with generator 
$-A$ then it follows from Proposition~\ref{pdton1502} that 
$0 \not\in \sigma_\ess(A)$ if and only if 
there exists an $\varepsilon > 0$ such that $\sigma(A) \cap [0,\varepsilon) = \{ 0 \} $.
In the next theorem we reformulate this  
by saying that $S_t$ converges in the operator norm
as $t \to \infty$.

\begin{thm} \label{tdton1.504}
Let $S$ be a self-adjoint irreducible Markov semigroup on $L_2(\Gamma)$ with 
generator $-A$.
The following are equivalent.
\begin{tabeleq}
\item \label{tdton1.504-1}
$0 \not\in \sigma_\ess(A)$.
\item \label{tdton1.504-2}
$\lim_{t \to \infty} S_t = P$ in $\cl(L_2(\Gamma))$.
\item \label{tdton1.504-3}
There exists an $\varepsilon > 0$ such that 
$\|S_t - P\|_{\cl(L_2(\Gamma))} \leq e^{- \varepsilon t}$
for all $t > 0$.
\end{tabeleq}
In that case one also has $\lim_{t \to \infty} S^{(p)}_t = P$ in $\cl(L_p(\Gamma))$
for all $p \in (1,\infty)$.
\end{thm}
\proof\
`\ref{tdton1.504-1}$\Rightarrow$\ref{tdton1.504-3}'.
We consider the situation in the proof of Proposition~\ref{pdton1501},
which was obtained via a unitary transformation.
Since $S$ is irreducible one has $\dim \ker A = 1$.
Then the hypothesis $0 \not\in \sigma_\ess(A)$ implies that 
there exists an 
$\varepsilon > 0$ such that $\sigma(A) \cap [0,\varepsilon) = \{ 0 \} $.
Then $m(y) \geq \varepsilon$ for a.e.\ $y \in Y \setminus Y_0$.
Thus 
\[
\|S_t - P\|_{\cl(L_2(Y))}
= \|S_t - P\|_{\cl(L_2(Y \setminus Y_0))}
= \|e^{-t m}\|_{L_\infty(Y \setminus Y_0)}
\leq e^{-\varepsilon t}
\]
for all $t > 0$.

`\ref{tdton1.504-3}$\Rightarrow$\ref{tdton1.504-2}' is trivial.

`\ref{tdton1.504-2}$\Rightarrow$\ref{tdton1.504-1}'.
The space $H_1 = (I - P)(L_2(\Gamma))$ is invariant under $S$
and $\lim_{t \to \infty} \|S_t\|_{\cl(H_1)} = 0$.
Since $S$ is self-adjoint, by the spectral theorem, this implies that 
there exists an $\varepsilon > 0$ such that 
$\|S_t\|_{\cl(H_1)} \leq e^{- \varepsilon t}$ for all $t > 0$.
Again by the spectral theorem this implies~\ref{tdton1.504-1}.

Finally we assume that \ref{tdton1.504-2} is valid.
Let $p \in (1,2)$.
Let $\theta \in (0,1)$ be such that $\frac{1}{p} = \frac{\theta}{1} + \frac{1-\theta}{2}$.
Then 
\[
\|S^{(p)}_t - P\|_{\cl(L_p(\Gamma))}
\leq \|S^{(1)}_t - P\|_{\cl(L_1(\Gamma))}^\theta \, \|S^{(2)}_t - P\|_{\cl(L_2(\Gamma))}^{1-\theta}
\leq 2^\theta \, \|S^{(2)}_t - P\|_{\cl(L_2(\Gamma))}^{1-\theta}
\]
for all $t > 0$ since $S^{(1)}$ is a contraction semigroup.
Therefore $\lim_{t \to \infty} S^{(p)}_t = P$ in $\cl(L_p(\Gamma))$.
The proof for $p \in (2,\infty)$ is similar, or follows by a duality argument.\hfill$\Box$

\ruimte

The harmonic oscillator on a weighted space (see \cite{Dav2} Theorem~4.3.6)
shows that the last assertion is not true, in general, for $p=1$
even if $A$ has compact resolvent.

\section{The Dirichlet-to-Neumann operator on arbitrary domains} \label{Sdton2.5}

In this section we will define the Dirichlet-to-Neumann operator $D_0$ on
$L_2(\Gamma)$ as a self-adjoint operator, and we will show that $-D_0$ generates
a Markov semigroup.

\begin{definition}\label{defi1.1}
Let $u\in H^1(\Omega)$ and $\varphi \in L_2(\Gamma)$.
We say that $\varphi$ is a \textbf{trace} of $u$ if there exist 
$u_1,u_2,\ldots\in H^1(\Omega) \cap C(\overline \Omega)$ 
such that $\lim\limits_{n\to \infty} u_n=u$
in $H^1(\Omega)$ and $\lim\limits_{n\to \infty} u|_\Gamma=\varphi$ in
$L_2(\Gamma)$.
\end{definition}

It is well possible that there are different elements of $L_2(\Gamma)$
such that they are both a trace of the same element of $H^1(\Omega)$ 
(see Section~\ref{Sdton3}).
Clearly if $u \in H^1(\Omega)$ has a trace then $u \in \widetilde H^1(\Omega)$.

Next we define the normal derivative $\frac{\partial u}{\partial \nu}$ by the
Green's formula  as follows (cf.\ \cite{AMPR} \cite{ArM} for 
the case that $\Omega$ has a Lipschitz boundary).
If $u\in L_{1,\loc}(\Omega)$, then we denote by $\Delta u \in \cd (\Omega)^\prime$
the distributional Laplacian applied to $u$.

\begin{definition}\label{defi1.2}
Let $u \in H^1(\Omega)$ be such that $\Delta u \in L_2(\Omega)$.
We say that $u$ has a normal derivative in $L_2(\Gamma)$ if there  exists a
$\psi \in L_2(\Gamma)$ such that
\begin{equation}\label{equ1.1}
\int_{\Omega} (\Delta u) \, v  + \int_{\Omega} \nabla u \cdot \nabla v
 = \int_\Gamma \psi \, v
\end{equation}
for all $v \in H^1(\Omega) \cap C(\overline \Omega)$.
In that case $\psi$ is
uniquely determined by (\ref{equ1.1}), we write 
$\frac{\partial u}{\partial \nu}:=\psi$ and call $\psi$ the 
\textbf{normal derivative} of $u$.
\end{definition}

To see uniqueness of the weak normal derivative, observe that by the
Stone--Weierstra{\ss} theorem the space $\{ v|_\Gamma : v\in \cd(\Ri^d)\}$ is
dense in $C(\Gamma)$ for the uniform norm and therefore also in $L_2(\Gamma)$.

\smallskip

Now  we are able to define the Dirichlet-to-Neumann operator $D_0$ on
$L_2(\Gamma)$.
It is part of the following theorem that the operator
$D_0$ is well defined, i.e.\ univocal, even though 
an $H^1(\Omega)$ function might have different functions
in $L_2(\Gamma)$ as a trace.

\begin{thm}\label{thm1.3}
There exists an operator $D_0$ on $L_2(\Gamma)$
such that the following holds.
Given $\varphi,\psi \in L_2(\Gamma)$ one has 
$\varphi \in D(D_0)$ and $D_0 \varphi = \psi$ if and only if there exists 
a $u \in H^1(\Omega)$
satisfying
\begin{itemize}
 \item 
$\Delta u = 0$,
\item 
$\varphi$ is a trace of $u$, and,
\item 
$u$ has a normal derivative in $L_2(\Gamma)$ and $\displaystyle\frac{\partial u}{\partial \nu} = \psi$.
\end{itemize}
Moreover, the operator $D_0$ is positive and self-adjoint.
\end{thm}

Here and in the sequel we always consider the operator $\Delta$ in the 
distributional sense.

For the proof of Theorem~\ref{thm1.3} we will need a generation theorem proved
recently in \cite{AE2} which is valid for arbitrary sectorial forms (without any
closability condition).
We recall a special case of it.

\begin{thm}  \label{thm1.4}
Let $D(a)$ be a real vector space and let $a\colon D(a) \times D(a) \to \Ri$ be
bilinear symmetric such that $a(u):=a(u,u) \geq 0$ for all $u\in D(a)$.
Let $H$ be a {\rm (}real\/{\rm )} Hilbert space and let $j\colon D(a)\to H$ be linear with dense
image.
Then there exists an operator $A$ on $H$ such that for all
$\varphi, \psi \in H$ one has $\varphi \in D(A)$ and $A\varphi = \psi$ if and
only if there exists a sequence $u_1,u_2,\ldots \in D(a)$ such that
\begin{tabel}
 \item
$\lim\limits_{n,m\to \infty} a(u_n-u_m)=0$,
\item
$\lim\limits_{n\to \infty} j(u_n)=\varphi$ in $H$, and,
\item
$\lim\limits_{n\to \infty} a(u_n,v)=(\psi, j(v))_H$ for all $v\in D(a)$.
\end{tabel}
Moreover, $A$ is positive and self-adjoint.
\end{thm}
\proof\
See \cite{AE2}, Theorem 3.2 and Remark~3.5.\hfill$\Box$

\ruimte

We call $A$ \textbf{the operator associated with} $(a,j)$.
Besides Theorem~\ref{thm1.4}, for the proof of Theorem~\ref{thm1.3},
we need the following remarkable inequality 
due to Maz'ya: There
exists a constant $c_M\geq 0$ such that
\begin{equation}\label{equ1.2}
\int_\Omega |u|^2 
\leq c_M \Big( \int_\Omega |\nabla u|^2 + \int_\Gamma |u|^2 \Big)
\end{equation}
for all $u \in H^1(\Omega) \cap C(\overline \Omega)$.
It was Daners (\cite{Daners2}) who showed how this inequality can
be used efficiently for elliptic and parabolic problems.
In fact, a stronger inequality is valid.
It follows from Example 3.6.2/1 and Theorem 3.6.3 in \cite{Maz} and (19) in \cite{AW2}
that there exists a constant $c_M' > 0$ such that 
\begin{equation} \label{equ1.3}
\Big( \int_\Omega |u|^q \Big)^{2/q}  
\leq c_M' \Big( \int_\Omega |\nabla u|^2 + \int_\Gamma |u|^2 \Big) 
\end{equation}
for all $u \in H^1(\Omega) \cap C(\overline \Omega)$, where
$q=\frac{2d}{d-1}$.
This inequality implies the following important compactness property 
(see \cite{Maz} Corollary 4.11.1/3).

\begin{prop}\label{prop1.5}
The space $H^1(\Omega)\cap C(\overline \Omega)$ with norm
\[
\|u\|^2
=\int_\Omega |\nabla u|^2 + \int_\Gamma |u|^2 
\]
is compactly embedded into $L_2(\Omega)$.
\end{prop}

In the proof of Theorem~\ref{thm1.3} we need the following form.
Define the form $\ell$ with form domain 
$D(\ell) = H^1(\Omega) \cap C(\overline \Omega)$ by 
\[
\ell(u,v) = \int_\Omega \nabla u\cdot \nabla v 
.  \]
The form $\ell$ is used throughout this paper.

\ruimte

\noindent
\textbf{Proof of Theorem~\ref{thm1.3}.}
Let $H=L_2(\Gamma)$.
Let $j\colon  D(\ell)\to L_2(\Gamma)$ be defined by $j(u)=u|_\Gamma$.
Then clearly $j$ has dense range.
Denote by $A$ the operator associated with $(\ell,j)$ in
the sense of Theorem~\ref{thm1.4}.
We shall show that $A$ has the properties of $D_0$.

Let $\varphi,\psi \in L_2(\Gamma)$.

Assume that $\varphi \in D(A)$ and $A\varphi = \psi$.
Then there exists a
sequence $u_1,u_2,\ldots \in D(\ell)$ such that 
$\lim\limits_{n,m\to \infty} \int_\Omega | \nabla(u_n-u_m)|^2=0$,
$\lim\limits_{n\to \infty} u_n|_\Gamma = \varphi$ in $L_2(\Gamma)$ and 
\begin{equation}\label{equ1.4}
\lim\limits_{n\to \infty} \int_\Omega \nabla u_n \cdot \nabla v  
= \int_\Gamma \psi \, v 
\end{equation}
for all $v\in D(\ell)$.
It follows from Maz'ya's inequality (\ref{equ1.2}) that $(u_n)_{n\in\Ni}$ is a
Cauchy sequence in $H^1(\Omega)$.
Let $u:=\lim\limits_{n\to \infty} u_n$ in
$H^1(\Omega)$.
Then $\varphi$ is a trace of $u$, by definition.
Moreover, by
(\ref{equ1.4}) we have
\[
\int_\Omega \nabla u \cdot \nabla v  = \int_\Gamma \psi \, v
\]
for all $v\in D(\ell)$.
Taking $v\in C^\infty_c(\Omega)$ we see that $\Delta u=0$.
Consequently,
\[
\int_\Omega (\Delta u) \,  v  + \int_\Omega \nabla u \cdot \nabla v  =
\int_\Gamma \psi \, v 
\]
for all $v \in D(\ell)$.
Therefore $u$ has a normal derivative in $L_2(\Gamma)$ and 
$\frac{\partial u}{\partial \nu} = \psi$ by Definition~\ref{defi1.2}.

Conversely, 
suppose there exists a $u\in H^1(\Omega)$ such that
$\Delta u=0$, the function $\varphi$ is a trace of $u$,
the function $u$ has a normal derivative in $L_2(\Gamma)$ and 
$\frac{\partial u}{\partial \nu}=\psi$.
Then there exist $u_1,u_2,\ldots \in D(\ell)$ such that 
$\lim\limits_{n\to \infty} u_n=u$ in $H^1(\Omega)$ and 
$\lim\limits_{n\to \infty} u_n|_\Gamma = \varphi$ in $L_2(\Gamma)$.
It follows that $\lim\limits_{n,m\to \infty} \ell(u_n-u_m)=0$ and, since $\Delta u = 0$, 
\begin{eqnarray*}
\lim\limits_{n\to \infty} \ell(u_n,v)
&=& \lim\limits_{n\to \infty}
\int_\Omega \nabla u_n \cdot \nabla v  
= \int_\Omega \nabla u \cdot \nabla v 
= \int_\Omega \nabla u \cdot \nabla v  + \int_\Omega (\Delta u) \, v 
= \int_\Gamma \psi \, v
\end{eqnarray*}
for all $v \in D(\ell)$ by the definition of $\frac{\partial u}{\partial \nu}$.
Hence $\varphi \in D(A)$ and $A\varphi=\psi$.

Therefore the operator with the properties of $D_0$ is well defined and equals $A$.
In particular $D_0$ is positive and self-adjoint.
This completes the proof of Theorem~\ref{thm1.3}.\hfill$\Box$

\ruimte

In the proof of Theorem~\ref{thm1.3} we also proved the following important fact,
which will be used later.

\begin{prop} \label{pdton240}
If $j \colon D(\ell) \to L_2(\Gamma)$ is defined by $j(u) = u|_\Gamma$,
then $D_0$ is the operator associated with $(\ell,j)$.
\end{prop}

We now show that the semigroup generated by $-D_0$ is markovian.

\begin{prop} \label{pdton340}
The $C_0$-semigroup $S$ on $L_2(\Gamma)$ generated by $-D_0$ is markovian,
i.e.\ $S_t \geq 0$ and $S_t \one_\Gamma = \one_\Gamma$ for all $t > 0$.
\end{prop}
\proof\
First we prove that $S$ is positive.
Let $L_2(\Gamma)_+ = \{ \varphi \in L_2(\Gamma) : \varphi \geq 0 \} $ be the 
positive cone in $L_2(\Gamma)$.
The orthogonal projection from $L_2(\Gamma)$ onto $L_2(\Gamma)_+$ is given by 
$\varphi \mapsto \varphi^+$.
Let $u \in H^1(\Omega) \cap C(\overline \Omega)$.
Then $u^+ \in D(\ell)$ and $j(u^+) = (j(u))^+$.
Moreover, 
$\ell(u^+, u - u^+) = - \ell(u^+, u^-) = - \int_\Omega \nabla (u^+) \cdot \nabla (u^-) = 0$
since $D_j (u^+) = \one_{[u > 0]} D_j u$ and 
$D_j (u^-) = - \one_{[u < 0]} D_j u$.
Hence $S$ is positive by Remark 3.12 in\cite{AE2}.

Since $\one_\Gamma \in D(D_0)$ and $D_0 \one_\Gamma = 0$ it follows 
that $S_t \one_\Gamma = \one_\Gamma$ for all $t > 0$.\hfill$\Box$

\section{Uniqueness of the trace and irreducibility} \label{Sdton3}

In general an element of $H^1(\Omega)$ might have more than one trace.
This happens if and only if the vector space 
$\{\varphi \in L_2(\Gamma):\varphi \mbox{ is a trace of } 0\}$
of degenerate traces is non-trivial.
By \cite{AE2} Lemma 4.14 there  
exists a Borel set $\Gamma_\sigma \subset \Gamma$ such that 
\[
\{ \varphi \in L_2(\Gamma):\varphi \mbox{ is a trace of } 0\}
=L_2(\Gamma\setminus \Gamma_\sigma) 
.  \]
Thus for all $\varphi \in L_2(\Gamma)$ one has 
$\varphi \in L_2(\Gamma\setminus\Gamma_{\sigma} )$ if and only if there exist 
$u_1,u_2,\ldots \in H^1(\Omega)\cap C(\overline \Omega)$ such that $\lim\limits_{n\to \infty}
\|u_n\|_{H^1(\Omega)}=0$ and $\lim\limits_{n\to \infty} u_n|_\Gamma=\varphi$ in $L_2(\Gamma)$.
We say that {\bf the trace on $\Omega$ is unique} if 
$\sigma(\Gamma \setminus \Gamma_\sigma)=0$; 
i.e.\ if $L_2(\Gamma \setminus \Gamma_\sigma)=\{0\}$.
This is equivalent with the fact that every element of $H^1(\Omega)$ has at most 
one trace.

Note that if $\sigma(\Gamma \setminus \Gamma_\sigma) > 0$ then the 
space $\Gamma \setminus \Gamma_\sigma$ is non-atomic since 
$d \geq 2$ (see \cite{Fre2}, Exercise~264Yg).
Hence $\dim L_2(\Gamma \setminus \Gamma_\sigma) = \infty$ if 
$\sigma(\Gamma \setminus \Gamma_\sigma) \neq 0$.

It follows from the definition of the operator
$D_0$ that $L_2(\Gamma \setminus \Gamma_\sigma) \subset \ker D_0$.
We next characterize $\ker D_0$.
In the proof we use that $\Omega$ is connected.

\begin{prop} \label{pdton213}
One has $\ker D_0 = \Ri \, \one_\Gamma + L_2(\Gamma \setminus \Gamma_\sigma)$.
Hence if $\sigma(\Gamma \setminus \Gamma_\sigma) = 0$, 
then $0 \in \sigma_{\rm p}(D_0)$ with multiplicity~$1$
and if $\sigma(\Gamma \setminus \Gamma_\sigma) > 0$, then 
$0 \in \sigma_{\rm p}(D_0)$ with infinite multiplicity.
\end{prop}
\proof\
Let $\varphi \in \ker D_0$.
By Theorem~\ref{thm1.3} there exists a $u \in H^1(\Omega)$
such that $\Delta u = 0$, $\varphi$ is a trace of $u$
and $0$ is the normal derivative of $u$.
Then $\int_\Omega \nabla u \cdot \nabla v = 0$
for all $v \in H^1(\Omega) \cap C(\overline \Omega)$.
Approximating $u$ by elements in $H^1(\Omega) \cap C(\overline \Omega)$
gives $\int_\Omega |\nabla u|^2 = 0$.
Since $\Omega$ is connected, one deduces that $u$ is constant.
So $\ker D_0 \subset \Ri \, \one + L_2(\Gamma \setminus \Gamma_\sigma)$.
The reverse inclusion is clear.\hfill$\Box$

\ruimte

\noindent
{\bf Proof of Theorem~\ref{tdton101}}\hspace{5pt}\ \
Theorem~\ref{tdton101} is a consequence of Theorem~\ref{tdton1.503} and 
Propositions~\ref{pdton1502} and \ref{pdton213}.\hfill$\Box$

\ruimte

If $\Omega$ is a Lipschitz domain, then 
$H^1(\Omega) = \widetilde H^1(\Omega)$ and 
there exists a $c > 0$ such that 
\[
\int_\Gamma |u|^2 \leq c\|u\|^2_{H^1(\Omega)}
\]
for all $u\in H^1(\Omega)\cap C(\overline \Omega)$.
This implies in particular that the trace on $\Omega$ is unique.
For general $\Omega$ 
it follows immediately from this result that the trace on $\Omega$ is unique
whenever there exists a Borel set $\Lambda \subset \Gamma$
with $\sigma(\Gamma \setminus \Lambda)= 0$
such that for each point $z\in \Lambda$ there exists an $r > 0$ such
that $B(z,r)\cap \Gamma$ is a Lipschitz graph with $B(z,r)\cap \Omega$ on one side.

There is another characterization for the uniqueness of the trace on $\Omega$
which involves the relative capacity on $\Omega$.
If $A \subset \Gamma$ is any set, then the 
{\bf relative capacity} of $A$ with respect to $\Omega$ is 
introduced in \cite{AW2} by
\begin{eqnarray*}
\capp_\Omega A
= \inf \{ \|u\|_{H^1(\Omega)}^2
& : & u \in \widetilde H^1(\Omega) 
    \mbox{ and there exists an open } V \subset \Ri^d 
    \mbox{ such } \\*
& & \mbox{that } A \subset V 
    \mbox{ and } u \geq 1 \mbox{ a.e.\ on }
       \Omega \cap V \}
.  
\end{eqnarray*}
Again another characterization is in terms of the Laplacian on $\Omega$
with Robin boundary conditions.
Define the form $a_R$ with domain $D(a_R) = H^1(\Omega) \cap C(\overline \Omega)$
by
\[
a_R(u,v) = \int_\Omega \nabla u \cdot \nabla v + \int_\Gamma u \, v
.  \]
Then $D(a_R)$ is a pre-Hilbert space with norm 
$\|u\|_{a_R}^2 = a_R(u) + \|u\|_{L_2(\Omega)}^2$.
Our second characterization of uniqueness of the trace is as follows.

\begin{prop} \label{pdton303}
The following conditions are equivalent.
\begin{tabeleq}
\item \label{pdton303-1}
The trace on $\Omega$ is unique.
\item \label{pdton303-3}
The form $a_R$ is closable.
\item \label{pdton303-2}
For every Borel set $B \subset \Gamma$ with 
$\capp_\Omega B = 0$ one has $\sigma(B) = 0$.
\end{tabeleq}
\end{prop}
\proof\
`\ref{pdton303-1}$\Rightarrow$\ref{pdton303-3}'.
Let $u_1,u_2,\ldots \in D(a_R)$ be a Cauchy sequence in $D(a_R)$ with 
$\lim u_n = 0$ in $L_2(\Omega)$.
Then $u_1,u_2,\ldots$ is a Cauchy sequence in $H^1(\Omega)$ and 
$u_1|_\Gamma, u_2|_\Gamma,\ldots$ is a Cauchy sequence in $L_2(\Gamma)$.
Hence $u := \lim u_n$ exists in $H^1(\Omega)$ and $\varphi := \lim u_n|_\Gamma$
exists in $L_2(\Gamma)$.
Then $u = 0$ since $\lim u_n = 0$ in $L_2(\Omega)$.
But the trace on $\Omega$ is unique.
So $\varphi = 0$ and consequently $\lim a_R(u_n) = 0$.
We have shown that $a_R$ is closable.

`\ref{pdton303-3}$\Rightarrow$\ref{pdton303-1}'.
Let $u_1,u_2,\ldots \in H^1(\Omega) \cap C(\overline \Omega)$, 
$\varphi \in L_2(\Gamma)$ and suppose that $\lim u_n = 0$ in $H^1(\Omega)$ and 
$\lim u_n|_\Gamma = \varphi$ in $L_2(\Gamma)$.
Then $u_1,u_2,\ldots$ is a Cauchy sequence in $D(a_R)$.
Moreover, $\lim u_n = 0$ in $L_2(\Omega)$ and $a_R$ is closable.
Therefore $\lim a_R(u_n) = 0$.
This implies that $\lim u_n|_\Gamma = 0$ in $L_2(\Gamma)$ and $\varphi = 0$.

`\ref{pdton303-3}$\Leftrightarrow$\ref{pdton303-2}'.
This is Theorem~3.3 in \cite{AW2}.\hfill$\Box$

\ruimte

The set $\Gamma_\sigma$ can also be described in a different way.
One says that $\sigma$ is {\bf admissible} if Property~\ref{pdton303-2}
of Proposition~\ref{pdton303} holds.
If $\sigma$ is not necessarily admissible, then there always exists a 
maximal admissible subset of $\Gamma$.
More precisely, the following is valid.

\begin{prop} \label{pdton302.2} 
There exists a Borel set $S \subset \Gamma$ such that 
\begin{tabel}
\item \label{pdton302.2-1} 
$\capp_\Omega(\Gamma \setminus S) = 0$ and 
\item \label{pdton302.2-2} 
if $B \subset \Gamma$ is a Borel set with 
$\capp_\Omega B = 0$, then 
$\sigma(B \cap S) = 0$.
\end{tabel}
\end{prop}
\proof\
See Proposition~3.6 in \cite{AW2}.\hfill$\Box$

\ruimte

It follows immediately from these two properties that the set 
$S$ in Proposition~\ref{pdton302.2} is $\sigma$-unique, i.e.\
if $S_1$ is another Borel set satisfying 
\ref{pdton302.2-1} and \ref{pdton302.2-2}, then $\sigma(S_1 \Delta S) = 0$.
If follows from the last paragraph of Section~3 in \cite{AW2}
that $\Gamma_\sigma$ equals $S$ up to $\sigma$-equivalence,
i.e.\ $\sigma(\Gamma_\sigma \Delta S) = 0$.

In \cite{AW2} Proposition 5.5
it is shown that always $\sigma(\Gamma_\sigma)>0$, without any regularity
assumption on the boundary (besides $\sigma(\Gamma)<\infty$).
Moreover, in \cite{AW2} Example 4.3 an example of a bounded connected open subset 
$\Omega \subset \Ri^3$ is given such that $\sigma(\Gamma) < \infty$ and 
$\sigma(\Gamma \setminus \Gamma_\sigma) > 0$.
A slightly easier example is as follows, which is a 
modification of an example at the end of Section~3 in \cite{BuG}.
It also has the property that $\widetilde H^1(\Omega) = H^1(\Omega)$.

\begin{voorb}[Uniqueness of the trace] \label{xdton450} \rm
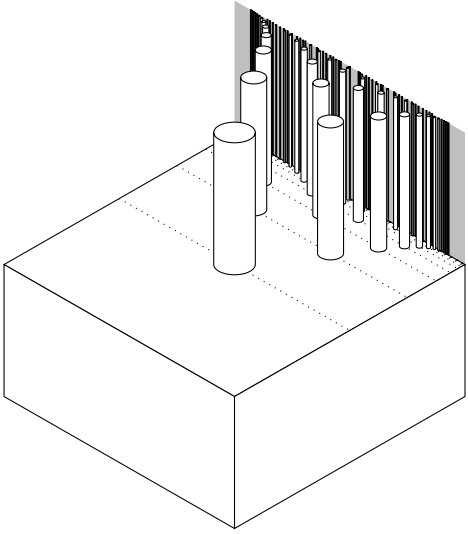
\begin{figure}[t]
\begin{minipage}[c]{0.48\textwidth}
\vspace{0pt}
\centering

\providecommand{\forrestdomainscale}{0.35}
\providecommand{\mycyl}[3]{
  \node[cylinder,draw,fill=white,shape aspect=0.5,rotate=90,minimum height={#2+#3/2},inner sep={#3/2},xshift={(#2-#3/2)/2}]
    at #1 {};
}
\begin{tikzpicture}[scale=\forrestdomainscale,every node/.style={scale=\forrestdomainscale}] 
    \def\tzrows{13}

\newcommand\xaxis{210}
\newcommand\yaxis{-30}
\newcommand\zaxis{90}

\coordinate(xone) at (\xaxis:1);
\coordinate(yone) at (\yaxis:1);
\coordinate(zone) at (\zaxis:1);

\newcommand\topside[3]{
  \draw[fill=white]
    (\zaxis:#3) -- ++(\xaxis:#1) -- ++(\yaxis:#2) -- ++({\xaxis+180}:#1) -- cycle;
}
\newcommand\leftside[3]{
  \draw[fill=white]
    (\xaxis:#1) -- ++(\yaxis:#2) -- ++(\zaxis:#3) -- ++({\yaxis+180}:#2) -- cycle;
}
\newcommand\rightside[3]{
  \draw[fill=white]
    (\yaxis:#2) -- ++(\zaxis:#3) -- ++(\xaxis:#1) -- ++({\zaxis+180}:#3) -- cycle;
}
\newcommand\cube[3]{
  \topside{#1}{#2}{#3} \leftside{#1}{#2}{#3} \rightside{#1}{#2}{#3}
}

\begin{scope}[shift={($(0,0)-5*(xone)$)}]
    \fill[fill=lightgray]
    (0,0) -- ++(\yaxis:10) -- ++(\zaxis:5) -- ++({\yaxis+180}:10) -- cycle;
\end{scope}

\begin{scope}[shift={($(0,0)-5*(xone)-5*(zone)$)}]
\cube{10}{10}{5}
\end{scope}

\begin{scope}
    \foreach \row in {\tzrows,...,1} {
      \pgfmathsetmacro\xc{(10/2.0*((1.0-pow(0.5,\row))/0.5-1.0))}
      \draw[style=dotted] ({\xaxis+180}:{\xc}) -- +(\yaxis:10); 
      \foreach \tree in {1,...,\row} {
        \pgfmathsetmacro\yc{(10/(\row+1.0)*\tree)}
        \pgfmathsetmacro\diam{(2.5cm/exp(\row*ln(1.6)))}
        \mycyl{($(0,0)-\xc*(xone)+\yc*(yone)$)}{5cm}{\diam};
      };
    };
\end{scope}
\end{tikzpicture}
\end{minipage}
\captionof{figure}{An example of a domain where $\sigma(\Gamma\setminus\Gamma_\sigma)>0$. 
In fact, the whole gray rectangle belongs to $\Gamma\setminus\Gamma_\sigma$.}
\label{fig:forrestdomain}
\end{figure}
For all $(x_0,y_0) \in [0,1] \times [0,1]$ and $r > 0$
let 
\[
C(x_0,y_0\,;r)
= \{ (x,y,z) \in \Ri^3 : |(x - x_0, y - y_0)| \leq r \mbox{ and } z \in [0,1] \}
\]
be the closed cylinder with axis parallel to the $z$-axis, radius $r$,
height $1$ and standing
on $(x_0,y_0,0)$.
Let 
\[
\Omega = \Int \Big( ([0,1] \times [0,1] \times [-1,0]) \cup 
   \bigcup_{n=1}^\infty \bigcup_{k=1}^{n-1} C(2^{-n}, \tfrac{k}{n} \,; 4^{-n}) \Big)
.  \]
(See Figure~\ref{fig:forrestdomain}.)
Then $\Omega$ is bounded, connected and $\sigma(\Gamma) < \infty$.
For all $m \in \Ni$ define $u_m \in H^1(\Omega) \cap C(\overline \Omega)$ by
\[
u_m(x,y,z) = (0 \vee (3^m \, z) \wedge 1) \, \one_{[0,2^{-m} + 4^{-m}]}(x)
.  \]
Then 
$\|u_m\|_{H^1(\Omega)}^2
\leq \sum_{n=m}^\infty \pi \, n \, (4^{-2n} + 3^{2m} \, 4^{-2n})$ 
for all $m \in \Ni$, so $\lim u_m = 0$ in $H^1(\Omega)$.
Since $0 \leq u_m \leq 1$ for all $m \in \Ni$ it follows from the 
Lebesgue domination convergence theorem that 
$\lim u_m|_\Gamma = \one_{ \{ 0 \} \times [0,1] \times [0,1]}$ in 
$L_2(\Gamma)$.
So $ \{ 0 \} \times [0,1] \times [0,1] \subset \Gamma \setminus \Gamma_\sigma$, 
up to $\sigma$-equivalence.

Let $u \in H^1(\Omega) \cap L_\infty(\Omega)$.
Define $v_m := u(\one - u_m)$ for all $m \in \Ni$.
Then $v_m$ has a support in a subdomain of $\Omega$ with a Lipschitz boundary.
So $v_m \in \widetilde H^1(\Omega)$.
Clearly $\sup_m \|u \, u_m\|_{H^1(\Omega)} < \infty$.
Therefore the sequence $v_1,v_2,\ldots$ has a weakly convergent subsequence in 
$\widetilde H^1(\Omega)$.
Moreover, $\lim u \, u_m = 0$ in $L_2(\Omega)$.
Hence $u \in \widetilde H^1(\Omega)$.
Since $H^1(\Omega) \cap L_\infty(\Omega)$ is dense in $H^1(\Omega)$
it follows that $H^1(\Omega) \subset \widetilde H^1(\Omega)$.
Thus $H^1(\Omega) = \widetilde H^1(\Omega)$.
\end{voorb}

In the above example the trace does not exist for each $u \in \widetilde H^1(\Omega)$.
We do not know whether universal existence of a trace implies its 
uniqueness.
More precisely, suppose that every element of $\widetilde H^1(\Omega)$ has a trace.
Does this imply that the trace on $\Omega$ is unique?

\section{Mapping properties of the trace} \label{Sdton4}

Let $H^1_\sigma(\Omega)$ be the set of all $u \in H^1(\Omega)$
for which there exists a $\varphi \in L_2(\Gamma)$ such that $\varphi$ is 
a trace of $u$.
Obviously, $H^1(\Omega) \cap C(\overline \Omega) \subset H^1_\sigma(\Omega)$.
It follows from the definition of the space $H^1_\sigma(\Omega)$ and the 
set $\Gamma_\sigma$ that there exists a
unique map 
\[
\Tr \colon H^1_\sigma(\Omega) \to L_2(\Gamma_\sigma)
\]
such that $\Tr u$ is a trace of $u$ for all $u \in H^1_\sigma(\Omega)$.
Then $\Tr u = u|_{\Gamma_\sigma}$ a.e.\ for all 
$u \in H^1(\Omega) \cap C(\overline \Omega)$.
Since $\Tr u$ is a trace of $u$ it follows from  the Maz'ya inequality 
(\ref{equ1.2}) that
\begin{equation}
\int_\Omega |u|^2 
\leq c_M \Big( \int_\Omega |\nabla u|^2 + \int_\Gamma |\Tr u|^2 \Big)  \\
= c_M \Big( \int_\Omega |\nabla u|^2 + \int_{\Gamma_\sigma} |\Tr u|^2  \Big)
\label{eSdton4;6}
\end{equation}
for all $u \in H^1_\sigma(\Omega)$, where we used that 
$(\Tr u)|_{\Gamma \setminus \Gamma_\sigma} = 0$.
Hence one can define the norm $\|\cdot\|_{H^1_\sigma(\Omega)}$ on 
$H^1_\sigma(\Omega)$ by 
\[
\|u\|_{H^1_\sigma(\Omega)}^2
= \int_\Omega |\nabla u|^2 + \int_{\Gamma_\sigma} |\Tr u|^2
.  \]
Obviously $\Tr \colon H^1_\sigma(\Omega) \to L_2(\Gamma_\sigma)$ is 
continuous.
We emphasize that in general the map 
$\Tr \colon (H^1_\sigma(\Omega), \|\cdot\|_{H^1(\Omega)}) \to L_2(\Gamma_\sigma)$
is not continuous.
A counter example is in \cite{Daners2} Remark~3.5(f).
It follows from (\ref{eSdton4;6}) that the norm
$\|\cdot\|_{H^1_\sigma(\Omega)}$ is equivalent to the norm
\[
u \mapsto \Big( \|u\|_{H^1(\Omega)}^2 + \|\Tr u\|_{L_2(\Gamma_\sigma)}^2
          \Big)^{1/2}
. \]
In particular $H^1_\sigma(\Omega)$ is a Hilbert space with inner product
\[
(u,v)_{H^1_\sigma(\Omega)} 
   = \int_\Omega \nabla u \cdot \nabla v + \int_\Gamma \Tr u \, \Tr v
\]
and $H^1_\sigma(\Omega)$ is continuously embedded in $L_2(\Omega)$.

The aim of this section is to study the map $\Tr$.
Before doing so, in the following remark, we show how the space
$H^1_\sigma(\Omega)$ can be used to give an alternative description of the
Dirichlet-to-Neumann operator.

\begin{remarkn}
The space $D(\ell)$ has the norm
\[
u \mapsto \Big( \int_\Omega |\nabla u|^2 + \int_\Gamma |u|^2 \Big)^{1/2}
.  \]
First we describe the completion of $D(\ell)$.
Define $\Phi \colon D(\ell) \to H^1_\sigma(\Omega) \oplus L_2(\Gamma \setminus \Gamma_\sigma)$
by $\Phi(u) = (u,u|_{\Gamma \setminus \Gamma_\sigma})$.
Then $\Phi$ is an isometry with dense range.
Therefore the space $H^1_\sigma(\Omega) \oplus L_2(\Gamma \setminus \Gamma_\sigma)$
is `the' completion of $D(\ell)$ and we identify $D(\ell)$ with 
$\Phi(D(\ell))$ in the natural manner.
Define the form $\tilde \ell$ with form domain
$D(\tilde \ell) = H^1_\sigma(\Omega) \oplus L_2(\Gamma \setminus \Gamma_\sigma)$
by 
\begin{equation}
\tilde \ell((u,\varphi), (v,\psi))
= \int_\Omega \nabla u \cdot \nabla v
\label{eSdton4;50}
\end{equation}
and define the map 
$\tilde j \colon H^1_\sigma(\Omega) \oplus L_2(\Gamma \setminus \Gamma_\sigma) \to L_2(\Gamma)$
by $\tilde j(u,\varphi) = \Tr u + \varphi$.
Then $\tilde \ell$ and $\tilde j$ are the continuous extensions
of $\ell$ and $j$, where $j \colon D(\ell) \to L_2(\Gamma)$ is 
defined by $j(u) = u|_\Gamma$.
Therefore $D_0$ is the operator associated with $(\tilde \ell,\tilde j)$
by \cite{AE2} Proposition~3.3.
Hence if $\varphi,\psi \in L_2(\Gamma)$, then $\varphi \in D(D_0)$ and
$D_0 \varphi = \psi$ if and only if there exists a 
$u \in H^1_\sigma(\Omega) \oplus L_2(\Gamma \setminus \Gamma_\sigma)$
such that $\tilde j(u) = \varphi$ and 
\begin{equation}
\tilde \ell(u,v) 
= (\psi,\tilde j(v))_{L_2(\Gamma)}
\label{eSdton4;5}
\end{equation}
for all $v \in H^1_\sigma(\Omega) \oplus L_2(\Gamma \setminus \Gamma_\sigma)$.
The latter follows from \cite{AE2} Theorem~2.1.
Then it follows immediately from (\ref{eSdton4;5}) that the range of $D_0$ is contained
in $L_2(\Gamma_\sigma)$.
\end{remarkn}


We will need the following apparently weaker description of the trace.

\begin{lemma} \label{ldton211}
Let $u \in L_2(\Omega)$ and $\varphi \in L_2(\Gamma)$.
Suppose there exist $u_1,u_2,\ldots \in H^1(\Omega) \cap C(\overline \Omega)$
such that $\lim u_n = u$ weakly in $L_2(\Omega)$,
$\lim u_n|_{\Gamma_\sigma} = \varphi|_{\Gamma_\sigma}$ weakly in $L_2(\Gamma_\sigma)$
and $\sup \|u_n\|_{H^1(\Omega)} < \infty$.
Then $u \in H^1_\sigma(\Omega)$ and $\varphi$ is a trace of $u$.
In particular, $\Tr u = \varphi \, \one_{\Gamma_\sigma}$.
\end{lemma}
\proof\
The sequence $u_1,u_2,\ldots$ is bounded in $H^1(\Omega)$ and the 
sequence $u_1|_{\Gamma_\sigma},u_2|_{\Gamma_\sigma},\ldots$
is bounded in $L_2(\Gamma_\sigma)$.
Therefore the sequence $u_1,u_2,\ldots$ is bounded in $H^1_\sigma(\Omega)$.
Since the unit ball is weakly compact it follows that,
after passing to a subsequence if necessary,
the sequence $u_1,u_2,\ldots$ is weakly convergent in $H^1_\sigma(\Omega)$.
So $u \in H^1_\sigma(\Omega)$.
Since the map $\Tr$ is bounded from $H^1_\sigma(\Omega)$ into 
$L_2(\Gamma_\sigma)$, it is also weakly continuous.
Hence $\Tr u = \lim \Tr u_n = \lim u_n|_{\Gamma_\sigma} = \varphi|_{\Gamma_\sigma}$
weakly in $L_2(\Gamma_\sigma)$.
So $\varphi \, \one_{\Gamma_\sigma} = \varphi|_{\Gamma_\sigma}$ is a trace of $u$.
Moreover, $\varphi \, \one_{\Gamma \setminus \Gamma_\sigma}$ is a trace of $0$.
Then $\varphi$ is a trace of $u$.\hfill$\Box$

\ruimte

We collect some algebraic properties of the trace.

\begin{prop} \label{pdton411.4}
\mbox{}
\begin{tabel}
\item \label{pdton411.4-1}
If $u \in H^1_\sigma(\Omega) \cap L_\infty(\Omega)$, then 
$\Tr u \in L_\infty(\Gamma_\sigma)$ and 
$\|\Tr u\|_\infty \leq \|u\|_\infty$.
Moreover, there exist $u_1,u_2,\ldots \in H^1(\Omega) \cap C(\overline \Omega)$ such that 
$\|u_n\|_\infty \leq \|u\|_\infty$ for all $n \in \Ni$,
$\lim u_n = u$ in $H^1(\Omega)$ and $\lim u_n|_\Gamma = \Tr u$ in $L_2(\Gamma)$.
\item \label{pdton411.4-2}
The space $H^1_\sigma(\Omega) \cap L_\infty(\Omega)$ is an algebra
and $\Tr (u \, v) = (\Tr u) \, (\Tr v)$ for all 
$u,v \in H^1_\sigma(\Omega) \cap L_\infty(\Omega)$.
\end{tabel}
\end{prop}
\proof\
`\ref{pdton411.4-1}'. 
There exist $u_1,u_2,\ldots \in H^1(\Omega) \cap C(\overline \Omega)$ such that 
$\lim u_n = u$ in $H^1(\Omega)$ and $\lim u_n|_\Gamma = \Tr u$ in $L_2(\Gamma)$.
For all $n \in \Ni$ set $v_n = (-M) \vee u_n \wedge M \in H^1(\Omega) \cap C(\overline \Omega)$,
where $M = \|u\|_\infty$.
Then $\lim v_n = (-M) \vee u \wedge M = u$ in $H^1(\Omega)$.
Moreover, $\lim v_n = (-M) \vee (\Tr u) \wedge M$ in $L_2(\Gamma)$.
So $(-M) \vee (\Tr u) \wedge M$ is a trace of $u$ and 
$\Tr u = \Big( (-M) \vee (\Tr u) \wedge M \Big) \, \one_{\Gamma_\sigma} 
       = (-M) \vee (\Tr u) \wedge M$.
Then $|\Tr u| \leq M$ a.e.
Note that $\|v_n\|_\infty \leq \|u\|_\infty$ for all $n \in \Ni$.

`\ref{pdton411.4-2}'. 
Let $u,v \in H^1_\sigma(\Omega) \cap L_\infty(\Omega)$.
By Statement~\ref{pdton411.4-1} 
there exist $u_1,u_2,\ldots,v_1,v_2,\ldots \in H^1(\Omega) \cap C(\overline \Omega)$
such that $\lim u_n = u$ in $H^1(\Omega)$, 
$\lim u_n|_\Gamma = \Tr u$ in $L_2(\Gamma)$,
$\lim v_n = v$ in $H^1(\Omega)$, 
$\lim v_n|_\Gamma = \Tr v$ in $L_2(\Gamma)$, and, moreover, 
$\|u_n\|_\infty \leq \|u\|_\infty$ and $\|v_n\|_\infty \leq \|v\|_\infty$
for all $n \in \Ni$.
Then $u_n \, v_n \in H^1(\Omega) \cap C(\overline \Omega)$ for all $n \in \Ni$
and 
\[
\|u_n \, v_n\|_{H^1(\Omega)}
\leq \|u_n\|_{H^1(\Omega)} \, \|v_n\|_\infty 
     + \|u_n\|_\infty \, \|v_n\|_{H^1(\Omega)}
\leq \|u_n\|_{H^1(\Omega)} \, \|v\|_\infty 
     + \|u\|_\infty \, \|v_n\|_{H^1(\Omega)}
\]
for all $n \in \Ni$.
So $\sup \|u_n \, v_n\|_{H^1(\Omega)} < \infty$.
Moreover, $\lim u_n \, v_n = u \, v$ in $L_2(\Omega)$ and 
$\lim (u_n \, v_n)|_\Gamma = (\Tr u) \, (\Tr v)$ in $L_2(\Gamma)$.
Therefore Lemma~\ref{ldton211} implies that $u \, v \in H^1_\sigma(\Omega)$ and 
$\Tr(u \, v) = (\Tr u) \, (\Tr v)$.\hfill$\Box$

\ruimte

The next lemma is a reformulation of Proposition~\ref{prop1.5}.

\begin{lemma} \label{ldton209}
The space $H^1_\sigma(\Omega)$ is compactly embedded 
in $L_2(\Omega)$.
\end{lemma}
\proof\
Let $B = \{ u \in H^1(\Omega) \cap C(\overline \Omega) : 
             \int_\Omega |\nabla u|^2 + \int_\Gamma |u|^2 \leq 2 \} $.
By Proposition~\ref{prop1.5} there exists a set $K \subset L_2(\Omega)$
which is compact in $L_2(\Omega)$ such that $B \subset K$.
Let $u \in H^1_\sigma(\Omega)$ and suppose that $\|u\|_{H^1_\sigma(\Omega)} \leq 1$.
There are $u_1,u_2,\ldots \in H^1(\Omega) \cap C(\overline \Omega)$ such that 
$\lim u_n = u$ in $H^1(\Omega)$ and $\lim u_n|_\Gamma = \Tr u$ in $L_2(\Gamma)$.
Then $u_n \in B \subset K$ for large $n$ and $\lim u_n = u$ in $L_2(\Omega)$.
So $u \in K$.\hfill$\Box$

\ruimte

Clearly $H^1_0(\Omega) \subset \{ u \in H^1_\sigma(\Omega) : \Tr u = 0 \} $.
If $\Omega$ is a Lipschitz domain, then the converse is valid (see \cite{Alt} Lemma~A~6.10).
We next give sufficient conditions for the converse inclusion, which 
allow $\Omega$ to have a cusp.

\begin{prop} \label{pdton220}
Suppose there exists a closed subset $K$ of $\Gamma$ such that 
$\capp_\Omega K = 0$ and for all $z \in \Gamma \setminus K$ there exists 
an $r > 0$ such that $B(z,r) \cap \Gamma$ is a Lipschitz graph
with $B(z,r) \cap \Omega$ on one side.
Then 
\[
\{ u \in H^1_\sigma(\Omega) : \Tr u = 0 \}
= H^1_0(\Omega)
.  \]
\end{prop}
\proof\
We may assume that $K \neq \emptyset$.
Let $u \in H^1_\sigma(\Omega)$ and suppose that $\Tr u = 0$.
We may assume that $u$ is bounded.

Let $\varepsilon > 0$.
We first prove that there exists a $\psi \in \widetilde H^1(\Omega)$
such that $0 \leq \psi \leq \one_\Omega$ a.e., 
$\|\psi\|_{H^1(\Omega)} \leq \varepsilon$ and 
$u (\one - \psi) \in H^1_0(\Omega)$.
Define the measure $\mu$ on the Borel $\sigma$-algebra of 
$\overline \Omega$ by $\mu(A) = |A \cap \Omega|$.
Define the form $h$ on $L_2(\overline \Omega,\mu)$
with form domain $D(h) = \widetilde H^1(\Omega)$ and 
$h(v,w) = \int_\Omega \nabla v \cdot \nabla w$.
Then $h$ is a regular Dirichlet form on $L_2(\overline \Omega,\mu)$
and $H^1(\Omega) \cap C(\overline \Omega)$ is a special standard core
for $h$ in the sense of \cite{FOT}.
Moreover, the relative capacity is just the capacity in \cite{FOT}
with respect to the Dirichlet form $h$ on $L_2(\overline \Omega,\mu)$.
For all $m \in \Ni$ let 
\[
K_m = \{ x \in \overline \Omega : d(x,K) \leq \tfrac{1}{m} \} 
.  \]
Then $K_m$ is compact, $K_1 \supset K_2 \supset \ldots$ and 
$\bigcap_{m=1}^\infty K_m = K$.
So by \cite{FOT} Theorem~2.1.1 there exists an $m \in \Ni$ such that 
$\capp_\Omega K_m < \varepsilon$.
Next, by \cite{FOT} Lemma~2.2.7(ii) there exists a 
$\psi \in H^1(\Omega) \cap C(\overline \Omega)$ such that 
$\one_{K_m} \leq \psi \leq \one$ and $\|\psi\|_{H^1(\Omega)}^2 \leq \varepsilon$.
It is an elementary exercise to see that there exists an open set 
$\Omega'$ in $\Ri^d$ with Lipschitz boundary
such that $\Omega \setminus K_m \subset \Omega' \subset \Omega$.
Let $\Gamma' = \partial(\Omega')$.
If $x \in \overline{\Omega'}$, then $x \in \overline \Omega$.
If $x \not\in \partial \Omega \cup K_m$ then 
$x \in \Omega \setminus K_m \subset \Omega'$.
So $\Gamma' \subset \Gamma \cup K_m$.
By Proposition~\ref{pdton411.4}\ref{pdton411.4-1} there exist
$u_1,u_2,\ldots \in H^1(\Omega) \cap C(\overline \Omega)$ such that 
$\|u_n\|_\infty \leq \|u\|_\infty$ for all $n \in \Ni$,
$\lim u_n = u$ in $H^1(\Omega)$ and $\lim u_n|_\Gamma = \Tr u = 0$ in $L_2(\Gamma)$.
For all $n \in \Ni$ define 
$v_n = (u_n (\one - \psi))|_{\overline{\Omega'}} \in H^1(\Omega') \cap C(\overline{\Omega'})$
and define $v = (u (\one - \psi))|_{\overline{\Omega'}}$.
Then
\[
\int_{\Gamma'} |v_n|^2
= \int_{\Gamma'} |u_n (\one - \psi)|^2
\leq \int_\Gamma |u_n (\one - \psi)|^2 + \int_{K_m} |u_n (\one - \psi)|^2
\leq \int_\Gamma |u_n|^2
.  \]
So $\lim v_n|_{\Gamma'} = 0$ in $L_2(\Gamma')$.
Moreover, $\lim v_n = v$ in $L_2(\Omega')$ and 
$\sup \|v_n\|_{H^1(\Omega')} \leq \sup \|u_n (\one - \psi)\|_{H^1(\Omega)} < \infty$.
So by Lemma~\ref{ldton211} it follows that $\Tr_{\Omega'} v = 0$.
Since $\Omega'$ has a Lipschitz boundary it follows that 
$v \in H^1_0(\Omega') \subset H^1_0(\Omega)$.
Then $u (\one - \psi) \in H^1_0(\Omega)$.

Let $n \in \Ni$.
By the above there exists a $\psi_n \in \widetilde H^1(\Omega)$
such that $0 \leq \psi_n \leq 1$ a.e.,
$\|\psi_n\|_{H^1(\Omega)} \leq \frac{1}{n}$ and $u (\one - \psi_n) \in H^1_0(\Omega)$.
Then $\sup \|u(\one - \psi_n)\|_{H^1_0(\Omega)}
\leq \sup \|u(\one - \psi_n)\|_{H^1(\Omega)} < \infty$.
So $n \mapsto u(\one - \psi_n)$ has a weakly convergent subsequence 
in $H^1_0(\Omega)$.
Alternatively, 
\[
\|u \, \psi_n\|_2
\leq \|u\|_\infty \, \|\psi_n\|_{H^1(\Omega)}
\leq \tfrac{1}{n} \, \|u\|_\infty
\]
for all $n \in \Ni$, so $\lim u(\one - \psi_n) = u$ in $L_2(\Omega)$.
Therefore $u \in H^1_0(\Omega)$.\hfill$\Box$

\section{Existence of a trace on $\widetilde H^1(\Omega)$} \label{Sdton6}

Recall that the trace $\Tr$ is defined on the subspace $H^1_\sigma(\Omega)$
of $\widetilde H^1(\Omega)$ and that in general the norm on 
$H^1_\sigma(\Omega)$ is strictly larger than the norm induced from 
$\widetilde H^1(\Omega)$.
In this section we characterize whether a trace exists
on all of $\widetilde H^1(\Omega)$.

\ruimte

We say that {\bf $\Omega$ has property~(P)}
if there exists a $c > 0$ such that 
\[
\int_\Omega |u - \langle u \rangle_\Omega|^2
\leq c \int_\Omega |\nabla u|^2
\]
for all $u \in H^1(\Omega) \cap C(\overline \Omega)$, 
where $\langle u \rangle_\Omega = \frac{1}{|\Omega|} \int_\Omega u$ is the 
{\bf average} of $u$ on $\Omega$.

Let $\widehat{D_0}$ be the part of the operator $D_0$ in 
the space $L_2(\Gamma_\sigma)$.
Then $\widehat{D_0}$ is a positive self-adjoint operator on $L_2(\Gamma_\sigma)$.

\begin{thm} \label{tdton210}
The following conditions are equivalent.
\begin{tabeleq}
\item \label{tdton210-5}
$H^1_\sigma(\Omega) = \widetilde H^1(\Omega)$ as sets, i.e.\ every
element of $\widetilde H^1(\Omega)$ has a trace.
\item \label{tdton210-3}
There exists a $c > 0$ such that 
$\int_{\Gamma_\sigma} |u|^2 \leq c \int_\Omega |\nabla u|^2$
for all $u \in H^1(\Omega) \cap C(\overline \Omega)$ with 
$\int_{\Gamma_\sigma} u = 0$.
\item \label{tdton210-4}
There exists a $c > 0$ such that 
$\int_{\Gamma_\sigma} |u|^2 \leq c \Big( \int_\Omega |\nabla u|^2 + \int_\Omega |u|^2 \Big)$
for all $u \in H^1(\Omega) \cap C(\overline \Omega)$.
\item \label{tdton210-7}
$0 \not\in \sigma_{\ess}(\widehat{D_0})$.
\end{tabeleq}
Moreover, if one of these equivalent conditions holds,
the space $\widetilde H^1(\Omega)$ is 
compactly embedded in $L_2(\Omega)$ and $\Omega$ has property~{\rm (P)}.
\end{thm}
\proof\
`\ref{tdton210-5}$\Rightarrow$\ref{tdton210-4}'.
If \ref{tdton210-5} is valid, then the norms on the spaces 
$H^1_\sigma(\Omega)$ and $\widetilde H^1(\Omega)$ are equivalent by the 
closed graph theorem.
Since $\Tr$ is continuous on $H^1_\sigma(\Omega)$ this implies 
that \ref{tdton210-4} is valid.

`\ref{tdton210-4}$\Rightarrow$\ref{tdton210-5}'.
One always has $H^1_\sigma(\Omega) \subset \widetilde H^1(\Omega)$.
Therefore the implication follows from 
Lemma~\ref{ldton211}.

`\ref{tdton210-3}$\Rightarrow$\ref{tdton210-4}'.
Define $F \colon (D(\ell), \|\cdot\|_{H^1(\Omega)}) \to \Ri$ by
\[
F(u) = \int_{\Gamma_\sigma} u
.  \]
We first prove that $F$ is continuous.
In order to prove this, it suffices to show that $\ker F$ is closed
in $(D(\ell), \|\cdot\|_{H^1(\Omega)})$.
Let $u_1,u_2,\ldots \in \ker F$, $u \in D(\ell)$ and suppose that 
$\lim u_n = u$ in $H^1(\Omega)$.
Then for all $\varepsilon > 0$ there exists an $N \in \Ni$ such that 
$\int_\Omega |\nabla(u_n - u_m)|^2 \leq \varepsilon$ for all $n,m \geq N$.
If $c > 0$ is as in \ref{tdton210-3}, it follows that 
$\int_{\Gamma_\sigma} |u_n - u_m|^2 \leq c \, \varepsilon$ for all $n,m \geq N$, 
where we used that $\int_{\Gamma_\sigma} (u_n - u_m) = F(u_n - u_m) = 0$ 
for all $n,m \in \Ni$.
So the sequence $u_1|_{\Gamma_\sigma}, u_2|_{\Gamma_\sigma},\ldots$ is a Cauchy sequence in 
$L_2({\Gamma_\sigma})$.
Hence $u_1,u_2,\ldots$ is Cauchy sequence in $H^1_\sigma(\Omega)$.
Since the space $H^1_\sigma(\Omega)$ is a Hilbert space, the
Cauchy sequence converges.
Therefore there exists a $\tilde u \in H^1_\sigma(\Omega)$
such that $\lim u_n = \tilde u$ in $H^1_\sigma(\Omega)$.
Then $\lim u_n = \tilde u$ in $L_2(\Omega)$, so $u = \tilde u \in H^1_\sigma(\Omega)$.
But $\Tr$ is continuous on $H^1_\sigma(\Gamma)$.
So $\lim \Tr u_n = \Tr u$ in $L_2(\Gamma_\sigma)$.
Then 
\[
\int_{\Gamma_\sigma} u
= (\Tr u,\one_{\Gamma_\sigma})_{L_2(\Gamma_\sigma)} 
= \lim (\Tr u_n,\one_{\Gamma_\sigma})_{L_2(\Gamma_\sigma)} 
= \lim F(u_n) 
= 0
.  \]
So $\ker F$ is closed and $F$ is continuous.
Hence there is a $c' > 0$ such that 
$|\langle u|_{\Gamma_\sigma} \rangle_{\Gamma_\sigma}|^2 \leq c' \|u\|_{H^1(\Omega)}^2$
for all $u \in H^1(\Omega) \cap C(\overline \Omega)$, where 
$\langle \varphi \rangle_{\Gamma_\sigma} 
       = \frac{1}{\sigma(\Gamma_\sigma)} \int_{\Gamma_\sigma} \varphi$
denote the average of $\varphi$ for all $\varphi \in L_1(\Gamma_\sigma)$.
Finally, let $u \in H^1(\Omega) \cap C(\overline \Omega)$.
Then 
\begin{eqnarray*}
\int_{\Gamma_\sigma} |u|^2
& = & \int_{\Gamma_\sigma} |u - \langle u|_{\Gamma_\sigma} \rangle_{\Gamma_\sigma}|^2 
    + \int_{\Gamma_\sigma} |\langle u|_{\Gamma_\sigma} \rangle_{\Gamma_\sigma}|^2  \\
& \leq & c \int_\Omega |\nabla u|^2 
  + |\langle u|_{\Gamma_\sigma} \rangle_{\Gamma_\sigma}|^2 \, \sigma(\Gamma_\sigma)
\leq (c + c' \, \sigma(\Gamma_\sigma) ) 
      \Big( \int_\Omega |\nabla u|^2 + \int_\Omega |u|^2 \Big)
\end{eqnarray*}
and \ref{tdton210-4} is valid.

If \ref{tdton210-5} is valid, then $\widetilde H^1(\Omega)$ is 
compactly embedded in $L_2(\Omega)$ by Lemma~\ref{ldton209}.
Since $\Omega$ is connected, it follows that $\Omega$ has property~(P).

`{\rm \ref{tdton210-4}}$\Rightarrow${\rm \ref{tdton210-3}}'.
If $c > 0$ is as in \ref{tdton210-4}, then 
\[
\int_{\Gamma_\sigma} |u - \langle u|_{\Gamma_\sigma} \rangle_{\Gamma_\sigma}|^2
\leq \int_{\Gamma_\sigma} |u - \langle u \rangle_\Omega|^2
\leq c \Big( \int_\Omega |\nabla u|^2 + \int_\Omega |u - \langle u \rangle_\Omega|^2 \Big)
\]
for all $u \in H^1(\Omega) \cap C(\overline \Omega)$.
Since $\Omega$ has property~(P), the implication 
{\rm \ref{tdton210-4}}$\Rightarrow${\rm \ref{tdton210-3}} follows.

`\ref{tdton210-3}$\Rightarrow$\ref{tdton210-7}'.
Suppose  \ref{tdton210-7}
is not valid.
Then $0 \in \sigma_\ess(\widehat{D_0})$.
It follows from Proposition~\ref{pdton213} that for all $n \in \Ni$ 
there exists a $\varphi_n \in D(\widehat{D_0})$ such that 
$\int_{\Gamma_\sigma} \varphi_n = 0$ and 
$0 < (\widehat{D_0} \, \varphi_n, \varphi_n)_{L_2(\Gamma_\sigma)} 
    \leq \frac{1}{n} \int_{\Gamma_\sigma} |\varphi_n|^2$.
Next there exists a unique $u_n \in H^1_\sigma(\Omega)$ such that 
$\Tr u_n = \varphi_n$ and 
$\int_\Omega \nabla u_n \cdot \nabla v 
   = (D_0 \varphi_n, \Tr v)_{L_2(\Gamma)} 
   = (\widehat{D_0} \, \varphi_n, \Tr v)_{L_2(\Gamma_\sigma)}$
for all $v \in H^1_\sigma(\Omega)$.
Therefore 
$\int_\Omega |\nabla u_n|^2 = (\widehat{D_0} \, \varphi_n, \varphi_n)_{L_2(\Gamma_\sigma)} 
   \leq \frac{1}{n} \int_{\Gamma_\sigma} |\Tr u_n|^2$.
So \ref{tdton210-3} is not valid.
Therefore \ref{tdton210-3}$\Rightarrow$\ref{tdton210-7}.

`\ref{tdton210-7}$\Rightarrow$\ref{tdton210-3}'.
Let $\ell_c$ and $\widehat{\ell_c}$ be the closed positive symmetric forms associated with 
$D_0$ and $\widehat{D_0}$.
Since $0 \not\in \sigma_\ess(\widehat{D_0})$ it follows from Proposition~\ref{pdton213}
that there is a $\mu > 0$ such that 
$\widehat{\ell_c}(\varphi) \geq \mu \int_{\Gamma_\sigma} |\varphi|^2$
for all $\varphi \in D(\widehat{l_c})$ with $\int_{\Gamma_\sigma} \varphi = 0$.
So by \cite{AE2} Theorem~2.5 it follows that 
$\tilde \ell(u) = \ell_c(\Tr u) = \widehat{\ell_c}(\Tr u) \geq \mu \int_{\Gamma_\sigma} |\Tr u|^2$
for all $u \in (\ker \Tr)^\perp \subset H^1_\sigma(\Omega)$,
where the orthoplement is in the Hilbert space $H^1_\sigma(\Omega)$
and $\tilde \ell$ is as in (\ref{eSdton4;50}).
Now let $u \in H^1_\sigma(\Omega)$ with $\int_{\Gamma_\sigma} \Tr u = 0$.
Write $u = u_1 + u_2$ with $u_1 \in (\ker \Tr)^\perp$ and $u_2 \in \ker \Tr$.
Then $\int_{\Gamma_\sigma} \Tr u_1 = 0$.
Moreover, $\tilde \ell(u_1,u_2) = (u_1,u_2)_{H^1_\sigma(\Omega)} = 0$.
Therefore 
\[
\tilde \ell(u) 
= \tilde \ell(u_1) + \tilde \ell(u_2)
\geq \tilde \ell(u_1)
\geq \mu \int_{\Gamma_\sigma} |\Tr u_1|^2
= \mu \int_{\Gamma_\sigma} |\Tr u|^2
.  \]
So \ref{tdton210-3} is valid.
This completes the proof of the theorem.\hfill$\Box$

\ruimte

Theorem~\ref{tdton102} characterizes when every element 
of $\widetilde H^1(\Omega)$ has a unique trace.

\ruimte

\noindent
{\bf Proof of Theorem~\ref{tdton102}\hspace{5pt}}\ \ 
If \ref{tdton102-1} or \ref{tdton102-4} is valid then the semigroup $S$
is irreducible by Theorem~\ref{tdton101} and Proposition~\ref{pdton213}.
Therefore \ref{tdton102-1}$\Leftrightarrow$\ref{tdton102-4} follows from 
Theorem~\ref{tdton1.504}.

`{\rm \ref{tdton102-4}}$\Rightarrow${\rm \ref{tdton102-2.5}}'.
If \ref{tdton102-4} is valid, then the trace on $\Omega$ is unique.
Then \ref{tdton102-2.5} follows from Theorem~\ref{tdton210}.

`{\rm \ref{tdton102-2.5}}$\Rightarrow${\rm \ref{tdton102-2}}'.
This is similar to the proof of 
{\rm \ref{tdton210-3}}$\Rightarrow${\rm \ref{tdton210-4}}
in the proof of Theorem~\ref{tdton210}.

`{\rm \ref{tdton102-2}}$\Rightarrow${\rm \ref{tdton102-3}}'
is trivial.

`{\rm \ref{tdton102-3}}$\Rightarrow${\rm \ref{tdton102-4}}'.
If \ref{tdton102-3} is valid, then the trace on $\Omega$ is unique.
Then \ref{tdton102-4} follows from Theorem~\ref{tdton210}.\hfill$\Box$

\section{Compact trace} \label{Sdton7} 

In the previous section we investigated when the trace is bounded from $\widetilde H^1(\Omega)$
into $L_2(\Gamma)$.
Now we want to characterize when the trace is compact.

\begin{prop} \label{pdton206}
The following are equivalent.
\begin{tabeleq}
\item \label{pdton216-1}
The Dirichlet-to-Neumann operator $D_0$ has a compact resolvent.
\item \label{pdton216-2}
The map $j \colon D(\ell) \to L_2(\Gamma)$ defined by $j(u) = u|_\Gamma$
is compact.
\item \label{pdton216-3}
The trace on $\Omega$ is unique and the
map $\Tr$ is compact {\rm (}from $H^1_\sigma(\Omega)$ into $L_2(\Gamma)${\rm )}.
\item \label{pdton216-4}
Every element in $\widetilde H^1(\Omega)$ has a unique trace and the map 
$\Tr \colon \widetilde H^1(\Omega) \to L_2(\Gamma)$ is compact.
\end{tabeleq}
\end{prop}
\proof\
`\ref{pdton216-1}$\Leftrightarrow$\ref{pdton216-2}'.
Let $\ell_c$ be the closed positive symmetric form on $L_2(\Gamma)$ 
associated with $D_0$.
Then $D_0$ has compact resolvent if and only if the embedding
from $D(\ell_c)$ into $L_2(\Gamma)$ is compact.
Let $V$ be the completion of $D(\ell)$, where $D(\ell)$ has
the (usual) norm $u \mapsto \Big( \int_\Omega |\nabla u|^2 + \int_\Gamma |u|^2 \Big)^{1/2}$.
Let $\tilde j \colon V \to L_2(\Gamma)$ be the continuous extension of the map $j$.
Then $D(\ell_c) = \tilde j ( (\ker \tilde j)^\perp )$, 
with the quotient norm of $(\ker \tilde j)^\perp$ by 
\cite{AE2} Theorem~2.5.
Therefore the embedding from $D(\ell_c)$ into $L_2(\Gamma)$
is compact if and only if $\tilde j|_{(\ker \tilde j)^\perp} \colon (\ker \tilde j)^\perp \to L_2(\Gamma)$
is compact.
The latter map is compact if and only if $\tilde j$ is compact and clearly that is 
equivalent with the compactness of the map~$j$.

`\ref{pdton216-1}$\Rightarrow$\ref{pdton216-4}'.
If $D_0$ has compact resolvent then $0 \not\in\sigma_\ess(D_0)$.
Hence every element of $\widetilde H^1(\Omega)$ has a unique trace by 
Theorem~\ref{tdton102}.
Moreover, the norms on $\widetilde H^1(\Omega)$ and $H^1_\sigma(\Omega)$
are equivalent.
So by \ref{pdton216-2} the map 
$\Tr|_{H^1(\Omega) \cap C(\overline \Omega)} \colon 
    (H^1(\Omega) \cap C(\overline \Omega), \|\cdot\|_{\widetilde H^1(\Omega)})
   \to L_2(\Gamma)$ is compact.
Then \ref{pdton216-3} follows by density.

`\ref{pdton216-4}$\Rightarrow$\ref{pdton216-3}$\Rightarrow$\ref{pdton216-2}'
is trivial.\hfill$\Box$

\begin{cor} \label{cdton218}
If the Dirichlet-to-Neumann operator $D_0$ has compact resolvent, then 
$\widetilde H^1(\Omega)$ is compactly embedded in $L_2(\Omega)$.
\end{cor}

\section{Robin boundary conditions for the Laplacian} \label{Sdton8}

Finally we wish to consider Robin boundary conditions with a possibly 
negative measure.
For all $\beta \in \Ri$ define the symmetric densely defined form $a_\beta$
by 
\[
a_\beta(u,v) 
= \int_\Omega \nabla u \cdot \nabla v - \beta \int_\Gamma u \, v
\]
with form domain $D(a_\beta) = H^1(\Omega) \cap C(\overline \Omega)$.

\begin{prop} \label{pdton304}
\mbox{}
\begin{tabel}
\item \label{pdton304-1}
Every element of $\widetilde H^1(\Omega)$ has a unique trace
if and only if there exists a $\beta > 0$
such that the form $a_\beta$ is lower bounded.
\item \label{pdton304-2}
The map
$\Tr$ is compact if and only if for all $\beta > 0$
the form $a_\beta$ is lower bounded.
\end{tabel}
\end{prop}
\proof\
Statement~\ref{pdton304-1} is easy, by 
Theorem~\ref{tdton102}\ref{tdton102-2}$\Leftrightarrow$\ref{tdton102-3},
so it remains to prove
Statement~\ref{pdton304-2}.
`$\Rightarrow$'. 
This is as in \cite{ArM} Proposition~2.2.
`$\Leftarrow$'.
Let $u_1,u_2,\ldots \in H^1_\sigma(\Omega)$ and suppose that 
$\lim u_n = 0$ weakly in $H^1_\sigma(\Omega)$.
We shall show that $\lim \Tr u_n = 0$ in $L_2(\Gamma)$.
Let $\varepsilon > 0$.
There exists an $M \geq 0$ such that 
$\int_\Omega |\nabla u_n|^2 \leq M$ for all $n \in \Ni$.
Note that $\widetilde H^1(\Omega) = H^1_\sigma(\Omega)$, with equivalent
norms, by Statement~\ref{pdton304-1} and Theorem~\ref{tdton210}.
Choosing $\beta = \frac{M}{\varepsilon}$, it follows from the assumption 
that there exists a $c > 0$ such that 
\[
\int_\Gamma |\Tr u|^2 
\leq \frac{\varepsilon}{M} \int_\Omega |\nabla u|^2 + c \int_\Omega |u|^2
\]
first for all $u \in H^1(\Omega) \cap C(\overline \Omega)$ and then 
by continuity for all $u \in \widetilde H^1(\Omega)$.
Then 
\begin{equation}
\int_\Gamma |\Tr u_n|^2 
\leq \varepsilon + c \, \|u_n\|_{L_2(\Omega)}^2
\label{epdton304;1}
\end{equation}
for all $n \in \Ni$.
Since the embedding of the space $H^1_\sigma(\Omega)$ into $L_2(\Omega)$
is compact by Lemma~\ref{ldton209}, one deduces that $\lim u_n = 0$
strongly in $L_2(\Omega)$.
Therefore one deduces from (\ref{epdton304;1})  that 
$\limsup \|\Tr u_n\|_{L_2(\Gamma)}^2 \leq \varepsilon$
and the proposition follows.\hfill$\Box$

\ruimte

We suppose for the remaining part of this section 
that every element of 
$\widetilde H^1(\Omega)$ has a unique trace.
Let 
\[
\beta_0
= \sup \{ \beta > 0 : \mbox{the form } a_\beta \mbox{ is lower bounded} \}
\in (0,\infty]
.  \]
One has $\beta_0 = \infty$ if $\Omega$ is a Lipschitz domain, 
but in general $\beta_0 < \infty$,
see Example~\ref{xdton217}.
It follows that $a_\beta$ is lower bounded for all $\beta \in (-\infty,\beta_0)$.
Let $R^{(\beta)}$ be the associated operator.

\begin{prop} \label{pdton703}
Let $\beta \in (-\infty,\beta_0)$ and $u,f \in L_2(\Omega)$.
Then $u \in D(R^{(\beta)})$ and $R^{(\beta)} u = f$ if and only if 
$u \in \widetilde H^1(\Omega)$, $- \Delta u = f$,
$u$ has a normal derivative in $L_2(\Gamma)$
and $\frac{\partial u}{\partial \nu} = \beta \Tr u$.
\end{prop}
\proof\
`$\Rightarrow$'.
Let $\beta_1 \in (\beta,\beta_0)$.
Then $a_{\beta_1}$ is lower bounded, so there exists a $\gamma_1 > 0$ 
such that $a_{\beta_1}(u) \geq - \gamma_1 \, \|u\|_{L_2(\Omega)}^2$
for all $u \in H^1(\Omega) \cap C(\overline \Omega)$.
Then $\beta_1 \int_\Gamma |u|^2 \leq \int_\Omega |\nabla u|^2 + \gamma_1 \int_\Omega |u|^2$
and 
\begin{equation}
(\beta_1 - \beta) \int_\Gamma |u|^2 
\leq \int_\Omega |\nabla u|^2 - \beta \int_\Gamma |u|^2 + \gamma_1 \int_\Omega |u|^2
= a_\beta(u) + \gamma_1 \int_\Omega |u|^2
\label{epdton703;1}
\end{equation}
for all $u \in H^1(\Omega) \cap C(\overline \Omega)$.
There exists a Cauchy sequence $u_1,u_2,\ldots$ in $D(a_\beta)$ such that 
$\lim u_n = u$ in $L_2(\Omega)$ and $\lim a_\beta(u_n,v) = (f,v)_{L_2(\Omega)}$
for all $v \in H^1(\Omega) \cap C(\overline \Omega)$.
Then $\sup a_\beta(u_n) < \infty$ and $\sup \int_\Omega |u_n|^2 < \infty$.
So by (\ref{epdton703;1}) also $\sup \int_\Gamma |u_n|^2 < \infty$ and 
subsequently $\sup \int_\Omega |\nabla u_n|^2 < \infty$.
So $(u_n)_{n \in \Ni}$ is bounded in $H^1(\Omega)$ and 
$(u_n|_\Gamma)_{n \in \Ni}$ is bounded in $L_2(\Gamma)$.
Without loss of generality we may assume that the sequence 
$(u_n)_{n \in \Ni}$ is weakly convergent in $H^1(\Omega)$ and 
$(u_n|_\Gamma)_{n \in \Ni}$ is weakly convergent in $L_2(\Gamma)$.
Since $\lim u_n = u$ in $L_2(\Omega)$ it follows that $u \in H^1(\Omega)$.
Then $u \in H^1_\sigma(\Omega)$ by Lemma~\ref{ldton211}.
Therefore $u \in \widetilde H^1(\Omega)$ by Proposition~\ref{pdton304}\ref{pdton304-1}
and Theorem~\ref{tdton210}.
Then 
\[
\int_\Omega \nabla u \cdot \nabla v - \beta \int_\Gamma (\Tr u) \, v
= \lim_{n \to \infty} a_\beta(u_n,v)
= \int_\Omega f \, v
\]
for all $v \in H^1(\Omega) \cap C(\overline \Omega)$.
Therefore $- \Delta u = f$, $u$ 
has a normal derivative in $L_2(\Gamma)$ 
and $\frac{\partial u}{\partial \nu} = \beta \Tr u$.

`$\Leftarrow$'.
There exist $u_1,u_2,\ldots \in H^1(\Omega) \cap C(\overline \Omega)$ such that 
$\lim u_n = u$ in $H^1(\Omega)$ and $\lim u_n|_\Gamma = \Tr u$ in $L_2(\Gamma)$.
It follows from the definition of $\frac{\partial u}{\partial \nu}$ that 
\[
\lim_{n \to \infty} a_\beta(u_n,v)
= \int_\Omega \nabla u \cdot \nabla v - \beta \int_\Gamma (\Tr u) \, v
= \int_\Omega f \, v
\]
for all $v \in H^1(\Omega) \cap C(\overline \Omega)$.
Moreover, $\lim u_n = u$ in $L_2(\Omega)$ and $u_1,u_2,\ldots$ is a 
Cauchy sequence in $D(a_\beta)$.
So $u \in D(R^{(\beta)})$ and $R^{(\beta)} u = f$.\hfill$\Box$

\ruimte

If $\beta \in (-\infty,\beta_0)$ then $R^{(\beta)}$ has compact resolvent by 
Lemma~\ref{ldton209} and \cite{AE2} Lemma~2.7.
This has consequences for the Dirichlet-to-Neumann operator.

\begin{prop} \label{pdton305}
If $\beta \in (0,\beta_0)$, then $\dim \ker(D_0 - \beta I) < \infty$ and 
$\sigma_{\rm p}(D_0) \cap [0,\beta]$ is finite.
\end{prop}
\proof\
Let $N \in \Ni$, $\beta_1,\ldots,\beta_N \in (0,\beta]$ and 
$\varphi_1,\ldots,\varphi_N$ be an orthonormal system in $L_2(\Gamma)$
such that $D_0 \varphi_n = \beta_n \, \varphi_n$ for all 
$n \in \{ 1,\ldots,N \} $.
Then $\varphi_n \in L_2(\Gamma_\sigma)$ since $\beta_n \neq 0$.
For all $n \in \{ 1,\ldots,N \} $ let $u_n \in H^1_\sigma(\Omega)$
be the unique element such that $\Tr u_n = \varphi_n$ and 
\[
\int_\Omega \nabla u_n \cdot \nabla v
= \beta_n \int_\Gamma \varphi_n \, \Tr v
\]
for all $v \in H^1_\sigma(\Omega)$.
Then 
\[
\int_\Omega \nabla u_n \cdot \nabla u_m
= \beta_n \int_\Gamma \varphi_n \, \varphi_m
= \beta_n \, \delta_{nm}
\]
and 
\[
(u_n, u_m)_{H^1_\sigma(\Omega)} 
= (\beta_n + 1) \, \delta_{nm}
\]
for all $n,m \in \{ 1,\ldots,N \} $.
Therefore $u_1,\ldots,u_N$ is linearly independent in $H^1_\sigma(\Omega)$.
Let $\alpha_1,\ldots,\alpha_N \in \Ri$.
Then 
\[
a_\beta(\sum_{n=1}^N \alpha_n \, u_n)
= \sum_{n,m=1}^N \alpha_n \, \alpha_m 
    \Big( \int_\Omega \nabla u_n \cdot \nabla u_m - \beta \int_\Gamma \varphi_n \, \varphi_m \Big)
= \sum_{n=1}^N |\alpha_n|^2 (\beta_n - \beta)
\leq 0
.  \]
Therefore 
\[
\spann \{ u_1,\ldots,u_N \} 
\subset \{ u \in H^1_\sigma(\Omega) : a_\beta(u) \leq 0 \}
.  \]
Since $R^{(\beta)}$ has a compact resolvent, the right hand space is 
finite dimensional.
This proves the proposition.\hfill$\Box$

\section{Examples} \label{Sdton5}

In this section we give two striking examples of connected bounded 
open sets with a continuous boundary such that the 
Dirichlet-to-Neumann operator does not have compact resolvent.
In both examples the trace on $\Omega$ is unique.
In one example every element of $\widetilde H^1(\Omega)$ has a trace,
in the other one not.

We first collect some known properties for domains with continuous boundary.
The last two statements use that $\Omega$ is connected.

\begin{prop} \label{pdton208}
Suppose $\Omega$ has a continuous boundary.
Then one has the following.
\begin{tabel}
\item \label{pdton208-1}
The space $H^1(\Omega) \cap C^\infty(\overline \Omega)$ is 
dense in $H^1(\Omega)$.
So in particular $\widetilde H^1(\Omega) = H^1(\Omega)$.
\item \label{pdton208-2}
The space $H^1(\Omega)$ is compactly embedded in $L_2(\Omega)$.
\item \label{pdton208-3}
The space $\Omega$ has property~{\rm (P)}.
\end{tabel}
\end{prop}
\proof\
Statement~\ref{pdton208-1} is in \cite{Maz} Theorem~1.1.6/2
and Statement~\ref{pdton208-2} is in 
\cite{EdE} Theorem~V.4.17.
Therefore $\Omega$ has the (Neumann type) Poincar\'e property,
which is in this case property~(P).\hfill$\Box$

\ruimte

We do not know, though, whether the trace on $\Omega$ is unique
if $\Omega$ has continuous boundary.

\ruimte

In the first example we explicitly give an element of $\widetilde H^1(\Omega)$
which does not have a trace.

\begin{exam} \label{xdton501}
Let 
\[
\Omega = \{ (x,y) \in \Ri^2 : 0 < x < 1 \mbox{ and } -x^4 < y < x^4 \}
.  \]
Clearly the set $\Omega$ is open, connected and the $1$-dimensional Hausdorff measure
of the boundary of $\Omega$ is finite.
Also $\Omega$ has a continuous boundary.
Therefore $\widetilde H^1(\Omega) = H^1(\Omega)$ and 
$H^1(\Omega)$ is compactly embedded in $L_2(\Omega)$.
Moreover, the trace on $\Omega$ is unique since $\Gamma \setminus \{ (0,0) \} $ is locally
Lipschitz.
Define $u \colon \Omega \to \Ri$ by $u(x,y) = \frac{1}{x}$.
Then $u \in H^1(\Omega)$.
Since 
\[
\int_0^1 |u(x,x^4)|^2 \, \sqrt{1 + (4 x^3)^2} \, dx
= \infty
, \]
it follows that $u$ does not have a trace.
In particular the Dirichlet-to-Neumann operator does not have compact resolvent
by Proposition~\ref{pdton206}.

It follows that the semigroup $S$ generated by $-D_0$ is not compact.
Therefore $S_t$ does not have a bounded kernel for all $t > 0$.
Hence $S_t$ does not map $L_2(\Gamma)$ into $L_\infty(\Gamma)$.
Since $S$ is submarkovian, this implies that $S$ is not ultracontractive.
In particular, there does not exists a $q > 2$ such that $u \mapsto u|_\Gamma$ 
maps $H^1(\Omega) \cap C^\infty(\overline \Omega)$ into $L_q(\Gamma)$.
\end{exam}

The next estimate is used in Example~\ref{xdton217}, but is also of independent interest.

\begin{lemma} \label{ldton216}
Let $e_1,e_2 \in \Ri^2$ with $\|e_1\| = \|e_2\| = 1$ and $|(e_1,e_2)| \neq 1$.
Let $a,b > 0$ and set 
\[
\Omega = \{ s \, e_1 + t \, e_2 : s \in (0,a) \mbox{ and } t \in (0,b) \} 
.  \]
Then
\[
\int_0^a |u(s \, e_1)|^2 \, ds 
\leq \frac{1}{\sqrt{1 - |(e_1,e_2)|^2}} 
   \Big( \frac{2}{b} \int_\Omega |u|^2 + b \int_\Omega |\nabla u|^2 \Big)
\]
for all $u \in H^1(\Omega) \cap C(\overline \Omega)$.
\end{lemma}
\proof\
Let $s \in (0,a)$ and $t \in (0,b)$.
Then 
\[
u(s \, e_1)
= u(s \, e_1 + t \, e_2) - \int_0^t e_2 \cdot (\nabla u)(s \, e_1 + r \, e_2) \, dr
\]
and therefore 
\[
|u(s \, e_1)|^2
\leq 2 |u(s \, e_1 + t \, e_2)|^2 + 2 t \int_0^b |(\nabla u)(s \, e_1 + r \, e_2)|^2 \, dr
.  \]
Hence integrating with respect to $t$ over $(0,b)$ and dividing 
by $b$ yields
\[
|u(s \, e_1)|^2
\leq \frac{2}{b} \int_0^b |u(s \, e_1 + t \, e_2)|^2 \, dt
   + b \int_0^b |(\nabla u)(s \, e_1 + r \, e_2)|^2 \, dr
\]
and 
\begin{eqnarray*}
\int_0^a |u(s \, e_1)|^2 \, ds
& \leq & \frac{2}{b} \int_0^a \int_0^b |u(s \, e_1 + t \, e_2)|^2 \, dt \, ds
   + b \int_0^a \int_0^b |(\nabla u)(s \, e_1 + r \, e_2)|^2 \, dr \, ds  \\
& = & \frac{1}{\sqrt{1 - |(e_1,e_2)|^2}} 
   \Big( \frac{2}{b} \int_\Omega |u|^2 + b \int_\Omega |\nabla u|^2 \Big)
\end{eqnarray*}
by a change of variables.\hfill$\Box$

\begin{lemma} \label{ldton212}
Let $a \in (0,1]$.
Define 
\[
\Omega = 
\{ (x,y) \in \Ri^2 : 0 < y < a \mbox{ and } |x| < a^2 - a \, y \} 
.  \]
Let $V = \{ u \in H^1(\Omega) \cap C(\overline \Omega)
            : u|_{[-a^2,a^2] \times \{ 0 \} } = 0  \} $.
Then 
\[
\frac{1}{3}
\leq
\sup \Big\{ \frac{\|\Tr u\|_{L_2(\Gamma)}^2}{\|u\|_{H^1(\Omega)}^2} 
         : u \in V \setminus \{ 0 \} \Big\}
\leq 2
.  \]
\end{lemma}
\proof\
Define $u \colon \overline \Omega \to [0,\infty)$ by $u(x,y) = y$.
Then $u \in V$.
Moreover, 
$\int_\Omega |u|^2 = \frac{a^5}{6}$,
$\int_\Omega |\nabla u|^2 = a^3$ and 
$\int_\Gamma |\Tr u|^2 = \frac{2}{3} \, a^3 \, \sqrt{1+a^2}$.
Therefore $\|\Tr u\|_{L_2(\Gamma)}^2 \geq \frac{1}{3} \, \|u\|_{H^1(\Omega)}^2$.
This proves the first inequality.

Next let $u \in V$ and $t \in [0,a]$.
Then $u(a^2 - a \, t, t) = \int_0^t u_y(a^2 - a \, t, s) \, ds$.
So 
\[
|u(a^2 - a \, t, t)|^2 
\leq t \int_0^t |(\nabla u)(a^2 - a \, t, s)|^2 \, ds
\leq a \int_0^t |(\nabla u)(a^2 - a \, t, s)|^2 \, ds
.  \]
Hence 
\begin{eqnarray*}
\sqrt{1 + a^2} \int_0^a |u(a^2 - a \, t, t)|^2 \, dt
& \leq & a \, \sqrt{1 + a^2} \int_0^a \int_0^t |(\nabla u)(a^2 - a \, t, s)|^2 \, ds \, dt  \\
& = & \sqrt{1 + a^2} \int_{\Omega_+} |\nabla u|^2
,  
\end{eqnarray*}
where $\Omega_+ = \Omega \cap ((0,\infty) \times \Ri)$.
So $\int_\Gamma |\Tr u|^2 \leq \sqrt{1 + a^2} \int_\Omega |\nabla u|^2$, 
and the lemma follows.\hfill$\Box$

\ruimte

We next give an example of an open connected bounded set $\Omega$ in $\Ri^2$
with continuous boundary, such that every element of $H^1(\Omega)$
has a unique trace,
$H^1(\Omega)$ is compactly embedded in $L_2(\Omega)$,
but the Dirichlet-to-Neumann operator does not have compact resolvent.

\begin{exam} \label{xdton217}
Let $\Omega_0 = (-1,1) \times (-1,0)$
and for all $n \in \Ni$ let 
\[
\Omega_n 
= \{ (x,y) \in \Ri^2 : 0 < y < a_n \mbox{ and } |x - 2^{-n}| < a_n^2 - a_n \, y \} 
,  \]
where $a_n = 4^{-n}$.
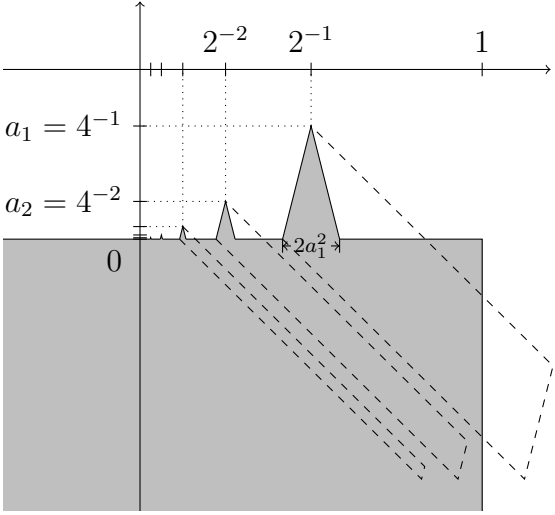
\begin{figure}[t]
\begin{minipage}[c]{0.48\textwidth}
\vspace{0pt}
\centering
\providecommand{\noncpttracescale}{4.5}
\begin{tikzpicture}[scale=\noncpttracescale] 
\begin{scope}
\clip (-0.4,-0.8) rectangle (1.25,0.7);

\draw[fill=lightgray] (-1,0)--(-1,-1)--(1,-1)--(1,0)--cycle;

\foreach \n in {1,...,5} {

\pgfmathsetmacro\yheight{pow(3,-\n)}
\pgfmathsetmacro\xmid{pow(2,-\n)}
\pgfmathsetmacro\xwidth{\yheight/4}
\pgfmathsetmacro\xleft{\xmid - \xwidth}
\pgfmathsetmacro\xright{\xmid + \xwidth}
\pgfmathsetmacro\xparbottom{\xleft + 1/sqrt(2)}
\pgfmathsetmacro\yparbottom{-1/sqrt(2)}
\pgfmathsetmacro\xparright{\xmid + 1/sqrt(2)}
\pgfmathsetmacro\yparright{\yheight - 1/sqrt(2)}

\draw[style={thick,color=lightgray}] (\xleft,0)--(\xright,0);
\draw[fill=lightgray] (\xleft,0)--(\xmid,\yheight)--(\xright,0);

\ifthenelse{\n<2}{
\draw (\xleft,0)--(\xleft,-0.04);
\draw (\xright,0)--(\xright,-0.04);
\path[<->,font=\scriptsize] (\xleft,-0.02) edge node[style={fill=lightgray,inner sep=1pt}] {$2a_1^2$} (\xright,-0.02);
}{}

\ifthenelse{\n<4}{
\draw[style=dashed] (\xleft,0)--(\xparbottom,\yparbottom)--(\xparright,\yparright)--(\xmid,\yheight);
\draw[style=dotted] (\xmid,\yheight)--(\xmid,0.5);
\draw[style=dotted] (\xmid,\yheight)--(0,\yheight);
}{}

\draw (\xmid,0.52) -- (\xmid,0.48);
\draw (-0.02,\yheight) -- (0.02,\yheight);

\ifthenelse{\n < 3}{
\draw (\xmid,0.52) node[anchor=south] {$2^{-\n}$};
\draw (-0.02,\yheight) node[anchor=east] {$a_{\n}=4^{-\n}$};
}{}
}

\draw[->] (-1.2,0.5) -- (1.2,0.5);
\draw[->] (0,-1.2) -- (0,0.7);

\foreach \x in {-1,1} {
\draw (\x,0.52) -- (\x,0.48);
\draw (\x,0.52) node[anchor=south] {$\x$};
}

\foreach \y in {-1,0} {
\draw (-0.02,\y) -- (0.02,\y);
\draw (-0.02,\y) node[anchor=north east] {$\y$};
}
\end{scope}
\end{tikzpicture}
\end{minipage}
\captionof{figure}{The domain in Example~\ref{xdton217}.}
\label{fig:noncpttrace}
\end{figure}
(See Figure~\ref{fig:noncpttrace}.)
Let $\displaystyle \Omega = \overline{\cup_{n=0}^\infty \Omega_n}\makebox[0pt]{\raisebox{4.5mm}[0pt]{\hspace{-14.5mm}$\scriptstyle\circ$}}$.
Then $\Omega$ is open, connected, the boundary is continuous
and $\sigma(\Gamma) < \infty$.

We show that every element of $\widetilde H^1(\Omega)$ has a unique trace
by showing that Condition \ref{tdton102-2} of Theorem~\ref{tdton102} is 
valid.
Since the set $\Omega_0$ is Lipschitz, it has the 
extension property. 
Therefore there exists a linear map 
$E \colon H^1(\Omega_0) \times C(\overline{\Omega_0}) \to H^1(\Ri^2) \times C(\Ri^2)$
and a constant $c_E > 0$ such that 
$(E u)|_{\Omega_0} = u$
and $\|E u\|_{H^1(\Ri^2)}^2 \leq c_E \, \|u\|_{H^1(\Omega_0)}^2$
for all $u \in H^1(\Omega_0) \times C(\overline{\Omega_0})$.

Let $u \in H^1(\Omega) \cap C(\overline \Omega)$.
Set $v = u|_{\overline{\Omega_0}}$ and $w = E v$.
Then $w \in H^1(\Ri^2) \times C(\Ri^2)$ and 
$(u-w)|_{\Omega_n} \in H^1(\Omega_n) \cap C(\overline{\Omega_n})$
with $(u-w)|_{[2^{-n} - a_n^2, 2^{-n} + a_n^2] \times \{ 0 \} } = 0$ for all $n \in \Ni$.
Then 
\begin{equation}
\int_\Gamma |u|^2
\leq \int_{\partial \Omega_0} |u|^2 + \int_{\Gamma \setminus \partial \Omega_0} |u|^2
\leq \int_{\partial \Omega_0} |u|^2 + 2 \int_{\Gamma \setminus \partial \Omega_0} |u-w|^2
   + 2 \int_{\Gamma \setminus \partial \Omega_0} |w|^2
.
\label{exdton217;1}
\end{equation}
We estimate the three terms in (\ref{exdton217;1}).

First, it follows from Lemma~\ref{ldton216} that 
\[
\int_{\partial \Omega_0} |u|^2
\leq 8 \|u\|_{H^1(\Omega_0)}^2
\leq 8 \|u\|_{H^1(\Omega)}^2
.  \]
Secondly, by Lemma~\ref{ldton212} one deduces that 
\begin{eqnarray*}
2 \int_{\Gamma \setminus \partial \Omega_0} |u-w|^2
& \leq & 2 \sum_{n=1}^\infty \int_{\partial \Omega_n} |u-w|^2
\leq 4 \sum_{n=1}^\infty \|u-w\|_{H^1(\Omega_n)}^2  \\
& \leq & 4 \|u-w\|_{H^1(\Omega)}^2
\leq 8 \|u\|_{H^1(\Omega)}^2 + 8 \|w\|_{H^1(\Omega)}^2
. 
\end{eqnarray*}
But 
\[
\|w\|_{H^1(\Omega)}^2
\leq \|w\|_{H^1(\Ri^2)}^2
\leq c_E \, \|v\|_{H^1(\Omega_0)}^2
\leq c_E \, \|u\|_{H^1(\Omega)}^2
.  \]
So 
\[
2 \int_{\Gamma \setminus \partial \Omega_0} |u-w|^2
\leq 8 (1 + c_E) \|u\|_{H^1(\Omega)}^2
.  \]
Therefore it remains to estimate the last term 
$\int_{\Gamma \setminus \partial \Omega_0} |w|^2$
in (\ref{exdton217;1}).
Let $n \in \Ni$ and set
\[
\Omega_n'
= \{ (2^{-n} - a_n^2, 0) 
   + s \, \tfrac{1}{\sqrt{1 + a_n^2}} (a_n,1)
   + t \, \tfrac{1}{\sqrt{2}} (1,-1) : s \in (0,a_n \, \sqrt{1+a_n^2}) \mbox{ and } t \in (0,1) \}
.  \]
Let $\Gamma_n^{(l)} = \partial \Omega_n \cap ( (-\infty,2^{-n}) \times (0,\infty) )$.
Then $\Gamma_n^{(l)} = \partial \Omega_n' \cap ( (-\infty,2^{-n}) \times (0,\infty) )$
and it follows from Lemma~\ref{ldton216} that 
\[
\int_{ \Gamma_n^{(l)} } |w|^2
\leq 4 \, \|w\|_{H^1(\Omega_n')}^2
\]
and therefore, again by disjointness,
\[
\sum_{n=1}^\infty \int_{ \Gamma_n^{(l)} } |w|^2
\leq 4 \, \|w\|_{H^1(\Ri^2)}^2
\leq 4 c_E \, \|u\|_{H^1(\Omega)}^2
.  \]
A similar estimate is valid on the right top boundary of $\partial \Omega_n$.
Combining these partial estimates with (\ref{exdton217;1}) one deduces 
that 
\[
\int_\Gamma |u|^2
\leq ( 16 + 24 c_E ) \|u\|_{H^1(\Omega)}^2
.  \]
So by Theorem~\ref{tdton102} every element of $\widetilde H^1(\Omega)$
has a unique trace.

For all $n \in \Ni$ define 
$u_n \in H^1(\Omega) \cap C(\overline \Omega)$ 
by $u_n(x,y) = y \, \one_{\overline{\Omega_n}}(x,y)$.
Then it follows from (the proof of) Lemma~\ref{ldton212} that
$\|\Tr u_n\|_{L_2(\Gamma)}^2 \geq \frac{1}{3} \, \|u_n\|_{H^1(\Omega)}^2$.
Since the norms on $H^1(\Omega)$ and $H^1_\sigma(\Omega)$ are 
equivalent and the functions $u_1,u_2,\ldots$ have disjoint support,
it follows that $\Tr$ is not compact.
Therefore the Dirichlet-to-Neumann operator does not have a compact resolvent.
\end{exam}

\subsection*{Acknowledgement}

The second named author is most grateful for the hospitality 
during a research stay at the University of Ulm.
Both authors wish to thank Manfred Sauter and Daniel Daners for their comments.
In addition we wish to thank Manfred Sauter for providing the two pictures.
Part of this work is supported by the Marsden Fund Council from Government funding, 
administered by the Royal Society of New Zealand.

\end{document}